\documentclass[english,10pt,onecolumn,draftcls]{IEEEtran}
\usepackage[T1]{fontenc}
\usepackage{geometry}
\usepackage[active]{srcltx}
\usepackage{color}
\usepackage{amsthm}
\usepackage{amsmath}
\usepackage{amssymb}
\usepackage{graphicx}
\usepackage{caption} 
\usepackage{subcaption}
\makeatletter
\theoremstyle{plain}
\newtheorem{thm}{\protect\theoremname}
\theoremstyle{plain}
\newtheorem{prop}[thm]{\protect\propositionname}
\theoremstyle{plain}
\newtheorem{lem}[thm]{\protect\lemmaname}

\@ifundefined{definecolor}
 {\usepackage{color}}{}
\usepackage{multicol}
\usepackage{flushend}
\usepackage[linesnumbered,ruled,boxed]{algorithm2e}
\usepackage{float}
\usepackage{footnote} 
\IEEEoverridecommandlockouts

\RequirePackage{colortbl, tabularx}
\@ifundefined{comment}{}
  {
   
  }%

\@ifundefined{showcaptionsetup}{}{%
 \PassOptionsToPackage{caption=false}{subfig}}
\usepackage{subfig}
\makeatother

\usepackage{babel}
\providecommand{\lemmaname}{Lemma}
\providecommand{\propositionname}{Proposition}
\providecommand{\theoremname}{Theorem}

\begin{document}
\global\long\def\Km{\zeta^{\max}}
\global\long\def\K{\mathcal{\zeta}}

\title{{\LARGE {\bf{Online Distributed ADMM on Networks:\\[-0.2in] Social Regret, Network Effect, and Condition Measures\footnote{A preliminary version of this work has appeared in the 2014 {\em IEEE Conference on Decision and Control}~\cite{Hosseini2014}.}}}}}

\author{Saghar Hosseini, Airlie Chapman, and Mehran Mesbahi\thanks{
The research of the authors was supported by the ONR grant N00014-12-1-1002 
and AFOSR grant FA9550-12-1-0203-DEF. The authors are with the Department of Aeronautics and Astronautics, University of Washington, WA 98105.  Emails: \{saghar, airliec, mesbahi\}@uw.edu.}}
\maketitle
\begin{abstract}
This paper examines online distributed Alternating Direction Method of Multipliers (ADMM). 
The goal is to distributively optimize a global objective function over a network of decision makers under 
linear constraints. 
The global objective function is composed of convex cost functions associated with each agent.
The local cost functions, on the other hand, are assumed to
have been decomposed into two distinct convex functions, one of which is revealed to the decision makers
over time and one known a priori. 
In addition, the agents must achieve consensus on the global variable that relates to the private local
variables via linear constraints.
In this work, we extend online ADMM to a distributed setting based on dual-averaging and
distributed gradient descent. We then propose a performance metric for such online distributed algorithms 
and explore the performance of the sequence of decisions generated by the algorithm as 
compared with the best fixed decision in hindsight.
This performance metric is called the social regret. A sub-linear upper bound on the 
social regret of the proposed algorithm
is then obtained that underscores the role of the underlying network topology and 
certain condition measures associated with the linear constraints. 
The online distributed ADMM algorithm
is then applied to a formation acquisition problem demonstrating
the application of the proposed setup in distributed robotics.

\end{abstract}


\section{Introduction}

Distributed convex optimization over networks arises in diverse application
domains, including multi-agent coordination, distributed estimation
in sensor networks, decentralized tracking, and event localization
\cite{Necoara2011756,Garcia2012}. A subclass of these problems can
be posed as optimization problems consisting of a composite convex
objective function subject to local linear constraints. This paper
examines two extensions of the well known Alternating Direction Method
of Multipliers (ADMM) algorithm~\cite{Lions1979} for solving this
class of problems. The first extension involves proposing two effective
means for distributed implementation of the ADMM algorithm. The second
extension pertains to addressing the situation where part of the cost
function has an online feature, representing uncertainties in the
cost incurred by each decision-maker prior to committing to a decision.
ADMM is an appealing approach that blends the benefits of augmented
Lagrangian and dual decomposition methods to solve the optimization
problem of the form, 
\begin{align}
\min_{x\in\chi,y\in Y}f(x)+\phi(y),\;\mbox{s.t. }Ax+By=c,\label{eq: ADMM}
\end{align}
where $f:\mathbb{R}^{d_{x}}\rightarrow\mathbb{R}$ and $\phi:{\mathbb{R}}^{d_{y}}\rightarrow{\mathbb{R}}$
are convex functions, and $\chi\subseteq{\mathbb{R}}^{d_{x}}$ and
$Y\subseteq{\mathbb{R}}^{d_{y}}$ are convex sets; $d_{x}$ and $d_{y}$
represent, respectively, the dimensions of the underlying Euclidean
spaces for the variables $x$ and $y$. %
ADMM has been extended to the scenarios where the cost function is
not known {\em a priori}. In other words, when the relevant decisions
are made, one part of the cost function might be varying with time,
or poorly characterized by a probability distribution, for example
due to uncertainties in the environment. In this case, the time varying
nature of this cost function is often signified by the notation $f_{t}$.
Such problem formulations fall under the class of online optimization
problems~\cite{Zinkevich2003}. Stochastic and online ADMM (O-ADMM)
have consequently been proposed to address this scenario in the context
of the following optimization problem at time $T>0$: 
\begin{align}
\min_{x\in\chi,y\in Y}\sum_{t=1}^{T}(f_{t}(x)+\phi(y)),\;\mbox{s.t. }Ax+By=c.\label{eq: online ADMM}
\end{align}
In this direction, stochastic ADMM has been introduced by Ouyang \textit{et
al.} \cite{Ouyang2012}, where an identical and independent distribution
for the uncertainties in the functions $f_{t}$ have been considered
and a convergence rate of $O(1/\sqrt{T})$ for convex functions has
been shown. The O-ADMM algorithms proposed in \cite{Wang2012a,Suzuki2013a}
also provide similar convergence rates without assumptions on the
distribution of uncertainties.

On the other hand, ADMM has been considered in the setting of distributed
convex optimization, particularly in the context of the consensus
problem~\cite{Boyd2010}, where agreement is required on each agent's
local variable $y_{i}$. In this case, the problem considered is of
the form, 
\begin{equation}
\min_{x\in\chi,y_{1},\dots,y_{n}\in Y}\sum_{i=1}^{n}\phi_{i}(y_{i}),\;\mbox{s.t. }x=y_{i},\,\mbox{ for }i=1,2,\ldots,n.\label{eq:consensus ADMM}
\end{equation}
In this consensus ADMM problem formulation, the local variables $y_{i}$'s
are required to reach consensus through the global variable $x$;
thus the linear constraint that ties the local variables to the global
variable is an equality. The consensus constraint set can also be
enforced through a network, where each agent coordinates on satisfying
the equality constraint with its neighboring agents. An important
distinction between the consensus ADMM and the problem of interest
in our work is that the objective functional in the former problem
setup does not explicitly have a term dictated by the global variable.
a natural extension of~\eqref{eq:consensus ADMM} for the solution
of distributed ADMM considered in this paper (by replacing the global
variable by its local copies and enforcing consensus) does not naturally
lead to a distributed solution strategy without resorting to a sequential
update~\cite{Wei2012} or inclusion of a fusion center~\cite{Boyd2010}.
works in distributed consensus ADMM include those based on stochastic
asynchronous edge based ADMM~\cite{Lutzeler2013, Wei2013b} and distributed
gradient descent~\cite{Yan2010,Duchi2012,Hosseini2013,Koppel2014},
where under the global objective \eqref{eq: online ADMM} and the
local objective \eqref{eq:consensus ADMM}, the rate of convergence
of $O(1/\sqrt{T})$ and $O(1/T)$ can be achieved, respectively. From
an algorithmic perspective, the approach proposed in this work is
also distinct from the stochastic asynchronous edge based ADMM proposed
in~\cite{Wei2013b,Lutzeler2013}. In particular, the embedding of
dual averaging in the distributed algorithm offers a privacy preserving
feature for the agents in the network. 
is, in the approach proposed in the present work, the local variables
remain private for each agent and only the dual variables are communicated
throughout the network. In applications such as cloud computing, the
privacy preserving feature of the proposed algorithm might be of great
interest for the security and reliability of the overall system.

Distributed ADMM has also been adopted for implementation on sensor
networks \cite{Schizas2008}. For example, Schizas \textit{et al.}
have proposed an algorithm that combines ADMM and block coordinate
descent that guarantees the sensors collectively converge to the maximum
likelihood estimate. The approach adopted by Schizas \textit{et al.}
is similar to one examined in~\cite{Boyd2010}, and as such, requires
an averaging step at each iteration and exchanging the primal variables
among the sensors.

ADMM has been examined in the context of optimization over certain
types of graphs. For example, Mota \textit{et al.} \cite{Mota2013}
have studied the ADMM consensus problem for connected bipartite graphs.
In particular, in \cite{Mota2013} it is shown that distributed ADMM
algorithm requires less communication between agents compared with
other algorithms for a given accuracy of the solution. Other works
in this area include that of Deng \textit{et al.} \cite{Deng2013}
which has proposed a proximal Jacobian ADMM suitable for parallel
computation. However, this method requires an all-to-all communication
over a complete graph in each iteration.

The main contribution of this work is twofold. First, we show that
both dual averaging and distributed gradient descent can seamlessly
be integrated in the ADMM setup, providing effective means for its
distributed implementation, or when the local variables are naturally
associated with decision-makers operating over a network. Second,
we show how network-level regret for such distributed ADMM can be
derived, highlighting the effect of the underlying network structure
on the performance of the algorithm and certain condition measures
for the linear constraints, when part of the cost structure has an
online character and is only revealed to the decision-makers over
time. As such, the paper extends and unifies some of the aforementioned
results on online and distributed ADMM. In the meantime, the paper
does not claim novelty in relation to developing a new class of ADMM
algorithms and instead builds on, and extends the existing ADMM iterations
for the purpose of its discussion. The paper considers the extension
of the optimization problem (\ref{eq: online ADMM}) of the form,
\begin{align}
\min_{\mbox{\scriptsize\ensuremath{\begin{array}{c}
x\in\mathcal{X},y_{1},\dots y_{n}\in Y\end{array}}}}\frac{1}{n}\sum_{t=1}^{T}\left[\sum_{i=1}^{n}\left\{ f_{i,t}(x)+\phi_{i}(y_{i})\right\} \right],\label{eq:online distributed ADMM}
\end{align}
\vspace{-0.1cm}
 
\[
\mbox{s.t. }A_{i}x+B_{i}y_{i}=c_{i}\mbox{ for }i=1,2,\ldots,n,
\]
involving a network of $n$ agents, each cooperatively solving for
the global optimal variable $x$ and the respective local variables
$y_{1},\dots y_{n}$. Here, the functions that compose problem \eqref{eq:online distributed ADMM}
are distributed, specifically only agent $i$ has access to functions
$f_{i,t}$, $\phi_{i}$, and its privately known local linear constraint.
of this problem include balancing sensing and communication in sensor
networks, analyzing large data sets in cloud computing, and cooperative
mission planning for a group of autonomous vehicles. The formulation
of cooperative forest firefighting using the optimization model (\ref{eq:online distributed ADMM})
and online distributed ADMM for its solution are discussed in \S\ref{sec:Example---Distributed}.



The outline of the paper is as follows. In \S\ref{sec:Background-and-Model},
the notation and a brief background on graphs and the regret framework
are presented. The optimization problem formulation and the network-level
measure of performance are introduced in \S\ref{sec:Problem-Statement}
followed by the description of the OD-ADMM algorithm and the corresponding
regret analysis in \S\ref{sec:Main-Result}. Then in \S\ref{sec:Example---Distributed},
the distributed formation acquisition problem is solved based on the
proposed algorithm, and simulation results are presented to support
the analysis. Finally, concluding remarks are provided in \S\ref{sec:Conclusion}.

\section{\label{sec:Background-and-Model}Background and Preliminaries}

In this section, we review basic concepts from graph theory and online
algorithms, as well as the relevant assumptions for our analysis.

The notation $v_{i}$ or $\left[v\right]_{i}$ denotes the $i$th
element of a column vector $v\in\mathbb{R}^{p}$. A unit vector $e_{i}$
denotes the column vector which contains all zero entries except $\left[e_{i}\right]_{i}=1$.
The vector of all ones will be denoted by $\mathbf{1}$. For a matrix
$M\in\mathbb{R}^{p\times q}$, $\left[M\right]_{ij}$ denotes the
element in its $i$th row and $j$th column. A doubly stochastic matrix
$P$ is a non-negative matrix with $\sum_{i=1}^{n}P_{ij}=\sum_{j=1}^{n}P_{ij}=1$.
For any positive integer $n$, the set $\{1,2,...,n\}$ is denoted
by $\left[n\right]$. The 2-norm,$1$-norm and infinity norm are denoted
by $||.||$, $||.||_{1}$, and $||.||_{\infty}$, respectively; the
dual norm of a vector $u$ in the normed space with the norm $\|.\|$
is defined as $||u||_{*}=\underset{||v||=1}{\sup}\left\langle u,v\right\rangle =||u||$,
where $\left\langle \cdot,\!\cdot\right\rangle $ denotes the underlying
inner product.

We denote the largest, second largest, and smallest singular values
of $Q\in\mathbb{R}^{n\times n}$ by $\sigma_{1}(Q)$, $\sigma_{2}(Q)$
and $\sigma_{n}(Q)$, respectively. A function $f:\chi\rightarrow\mathbb{R}$
is called $L$-Lipschitz continuous if there exists a positive constant
$L$ for which 
\begin{equation}
|f(u)-f(v)|\leq L\Vert u-v\Vert\;\;\mbox{for all}\; u,v\in\chi.\label{eq:L-Lipschitz}
\end{equation}
Although the dual of the 2-norm is the 2-norm itself, we derive some
of the bounds in our subsequent analysis using the notion of the dual
norm. The main reason is the connection between the Lipschitz continuity
of a function (in the native norm) and the boundedness of its subgradient
(by the Lipschitz constant) in the dual norm.


A graph is an abstraction for representing the interactions among
decision-makers, e.g., sensors and mobile robots. A weighted graph
$\mathcal{G}=\left(V,E,W\right)$ is defined by the node set $V$,
where the number of nodes in the graph is $|V|=n$. Nodes represent
the decision-makers in the network, and the edge set $E$ represents
the agents' interactions, that is, agent $i$ communicates with agent
$j$ if there is an edge from $i$ to $j$, i.e., $\left(i,j\right)\in E$.
In addition, a weight $\mbox{\ensuremath{w_{ji}\in W}}$ can be associated
with every edge $\left(i,j\right)\in E$ through the function $W:E\rightarrow\mathbb{R}$.
The neighborhood set of node $i$ is defined as $N(i)=\{j\in V|(i,j)\in E\}$.
One way to represent $\mathcal{G}$ is through the adjacency matrix
$A(\mathcal{G})$ where $\left[A(\mathcal{G})\right]_{ji}=w_{ji}$
for $\left(i,j\right)\in E$ and $\left[A(\mathcal{G})\right]_{ji}=0$,
otherwise. For a graph $\mathcal{G}$, $d_{i}$ is the weighted in-degree
of $i$ defined as $d_{i}=\sum_{\{j|(j,i)\in E\}}w_{ij}$. Another
matrix representation of $\mathcal{G}$ is the weighted graph Laplacian
defined as $L(\mathcal{G})=\Delta(\mathcal{G})-A(\mathcal{G})$, where
$\Delta(\mathcal{G})$ is the diagonal matrix of node in-degree's
$d_{i}$. If there exists a directed path between every pair of distinct
vertices, the graph $\mathcal{G}$ is referred to as strongly connected.
In this work, we assume that the inter agent communication between
the agents constitute a strongly connected graph, ensuring information
flow amongst the agents.


In online optimization, an online algorithm generates a sequence of
decisions $x_{t}$. At iteration $t$, the convex cost function $l_{t}$
remains unknown prior to committing to $x_{t}$. The feedback available
to the algorithm is the loss $l_{t}(x_{t})$ and its gradient. We
capture the performance of online algorithms by a standard measure
called regret. Regret measures how competitive the algorithm is with
respect to the best fixed solution. This best fixed decision, denoted
as $x^{*}$, is chosen with the benefit of hindsight. Formally, the
regret is defined as the difference between the incurred cost $l_{t}(x_{t})$
and the cost of the best fixed decision $l_{t}(x^{*})$ after $T$
iterations, i.e., 
\begin{equation}
R_{T}=\sum_{t=1}^{T}\left\{ l_{t}(x_{t})-l_{t}(x^{*})\right\} .\label{eq:regret_general}
\end{equation}

An online algorithm performs well if its regret grows sub-linearly
with respect to the number of iterations, i.e., 
\[
\lim_{T\rightarrow\infty}R_{T}/T=0.
\]
This implies that the average loss of the algorithm tends to the average
loss of the best fixed strategy in hindsight independent of the uncertainties
associated with the global cost.%
\footnote{The notion of regret is often received with a degree of skepticism
upon initial encounter. The basic idea is that if there is a positive
lower bound between the cost incurred by the algorithm and the best
fixed decision in hindsight, then the regret will grow linearly. A
sublinear regret implies that the algorithm has learned to match the
performance of the best fixed decision in hindsight.%
} We refer to \cite{Shalev-Shwartz2011b,SebastienBubeck2011,Hazan2011,Hazan2007a}
for further discussions on online algorithms and their regret analysis.

\section{Problem Statement\label{sec:Problem-Statement}}

In this section, we consider a large scale network of agents cooperatively
optimizing a global objective function. Let the communication geometry
amongst the $n$ decision-makers or agents, be denoted by the graph
$\mathcal{G}=(V,E)$. Each node $i\in V$ is an agent that communicates
with its neighbor $j\in N(i)$ through the edge $(i,j)\in E$. An
equivalent online distributed convex optimization problem to \eqref{eq:online distributed ADMM}
is as follows, 
\begin{align}
\min_{\begin{array}{c}
x\in\chi,y_{1},\dots,y_{n}\in Y\end{array}}\hspace{-0in}\sum_{t=1}^{T}F_{t}(x,y) & :=\sum_{t=1}^{T}\{f_{t}(x)+\frac{1}{n}\sum_{i=1}^{n}\phi_{i}(y_{i})\}\label{eq:objective fun}
\end{align}
subject to 
\begin{equation}
r_{i}(x,y_{i}):=A_{i}x+B_{i}y_{i}-c_{i}=0\;\mbox{for all }i\in[n],\label{eq:local linear constraints}
\end{equation}
where $f_{t}(x)=\frac{1}{n}\sum_{i=1}^{n}f_{i,t}(x)$, and $f_{i,t}:\mathbb{\mathbb{R}}^{d_{x}}\rightarrow\mathbb{R}$
and $\phi_{i}:\mathbb{\mathbb{R}}^{d_{y}}\rightarrow\mathbb{R}$ are
convex for each $i$. The matrices in the local linear constraints
are denoted as $A_{i}\in\mathbb{\mathbb{R}}^{m_{i}\times d_{x}}$,
$B_{i}\in\mathbb{\mathbb{R}}^{m_{i}\times d_{y}}$, and $c_{i}\in\mathbb{\mathbb{R}}^{m_{i}}$
at node $i\in[n]$ . We assume that $B_{i}^{T}$ is left invertible,
i.e., $\sigma_{m_{i}}(B_{i}B_{i}^{T})$ is non-zero, for all $i\in[n]$.
The functions $f_{i,t}$ and $\phi_{i}$ are further assumed to be
Lipschitz continuous with Lipschitz constants $L_{f}$ and $L_{\phi},$
respectively, that is, 
\begin{align*}
|f_{i,t}(u)-f_{i,t}(v)| & \leq L_{f}\Vert u-v\Vert\quad\mbox{for all }u,v\in\chi,\\
|\phi_{i}(u')-\phi_{i}(v')| & \leq L_{\phi}\Vert u'-v'\Vert\quad\mbox{for all }u',v'\in Y.
\end{align*}
The distributed nature of the optimization is illustrated in Figure
\ref{fig:Distributed-ADMM-problem}. 
assume that the Slater condition holds, namely that there exists $(x,y_{1},\cdots,y_{n})$
in the interior of $\chi\times Y\ldots\times Y$ 
that \eqref{eq:local linear constraints} is satisfied. This assumption
is naturally used in the analysis of the duality gap for deriving
bounds on the social regret. Moreover, we assume that the set of optimal
solutions of \eqref{eq:objective fun} 
nonempty and the finite optimum value is $\mathcal{P}^{*}$. 
diameter of the set $\chi$, defined as $\mbox{{\bf diam}}(\chi)=\sup_{x,x'\in\chi}||x-x'||$,
is assumed to be finite and denoted by $D_{\chi}$. 

\begin{center}
\begin{figure}[tb]
\begin{centering}
\includegraphics[width=0.4\columnwidth]{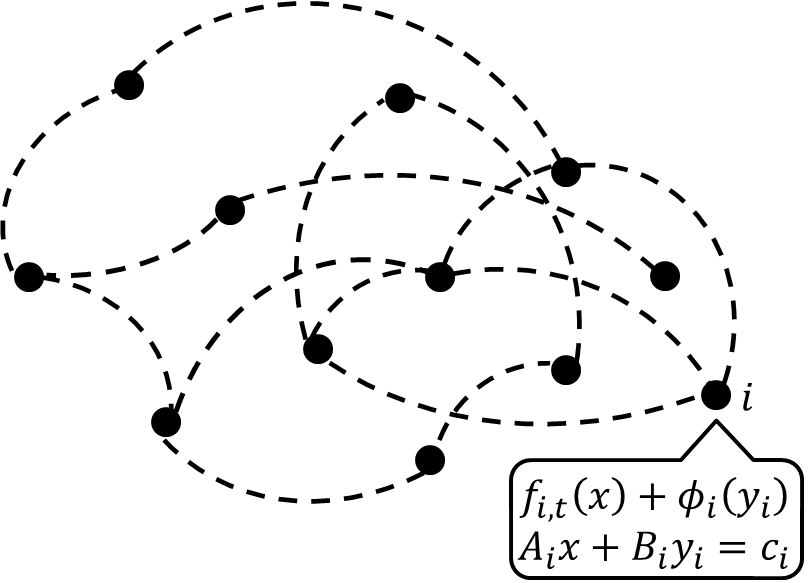} 
\par\end{centering}

\protect\protect\caption{Distributed ADMM problem over a network; each agent operates based
of the local objective $f_{i,t}(x)+\phi_{i}(y_{i})$ and the local
linear constraint $A_{i}x+B_{i}y_{i}=c_{i}$.\label{fig:Distributed-ADMM-problem}}
\end{figure}

\par\end{center}

The \textit{local} decisions made by agent $i$ is represented by
the optimization variables $x_{i}\in\mathcal{X}$ and $y_{i}\in Y$;
note that we allow the agents to have a local (not necessary exact)
version of the global variable $x$, namely $x_{i}$. In addition,
we assume that subgradients $\partial f_{i,t}(x)$ can be computed
for every $x\in\chi$. In the online setting, based on the available
local information, each decision maker $i$ selects a global variable
$x_{i,t}\in\chi$ and local variable $y_{i,t}\in Y$, at time $t$.
The cost $f_{i,t}(x_{i,t})$ is then revealed to this agent after
its local decision $x_{i,t}$ has been committed to at time $t$.

\subsection{Regret for Constrained Optimization }

We now examine a measure for evaluating the performance of OD-ADMM
based on variational inequalities. 
measure is inspired by the convergence analysis of Douglas-Rachford
ADMM presented in \cite{He2012}.

\begin{figure}[tb]
\centering{}\includegraphics[width=0.5\textwidth]{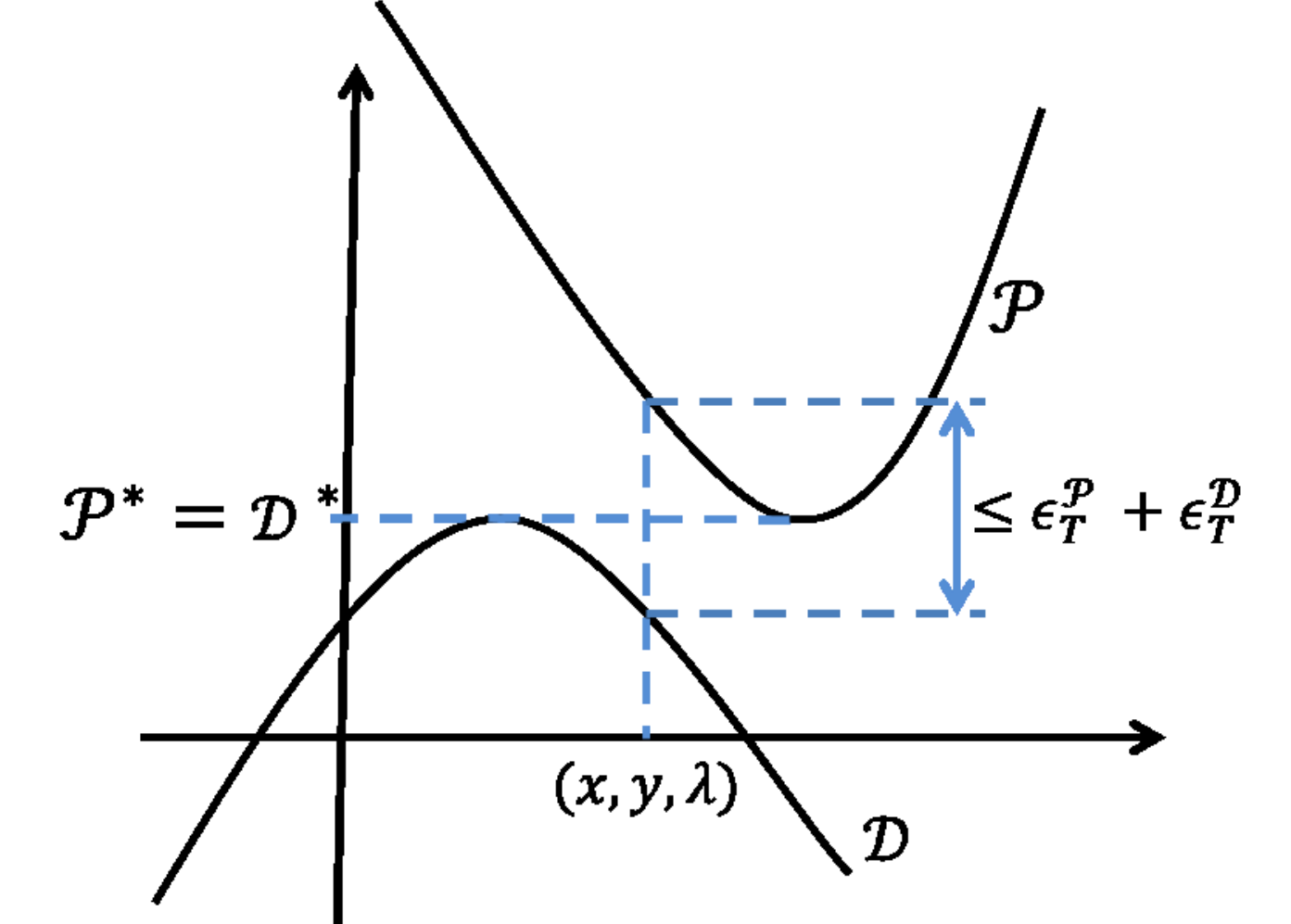}\protect
\protect\caption{\label{fig:duality gap}This figure illustrates the duality gap at
$w=\left(x,y,\lambda\right)\in\Omega$.}
\end{figure}

Consider the Lagrangian for the constrained optimization problem \eqref{eq:objective fun}
as 
\begin{equation}
\mathcal{L}_{T}(x,y,\lambda)=\sum_{t=1}^{T}\left\{ f_{t}(x)+\frac{1}{n}\sum_{i=1}^{n}(\phi_{i}(y_{i})+\left\langle \lambda_{i},r_{i}(x,y_{i})\right\rangle )\right\} ,\label{eq: Lagrangian}
\end{equation}
where $x\in\mathcal{X}$ and $y_{i}\in Y$, as well as assuming $\lambda_{i}\in\mathbb{\mathbb{R}}^{m_{i}}$,
for all $i\in\left[n\right]$. Then, the Lagrange dual function is
defined as 
\begin{equation}
\mathcal{D}(\lambda)=\inf_{x\in\chi,y_{i}\in Y}\mathcal{L}_{T}(x,y,\lambda),\label{eq:dual func}
\end{equation}
implying that $\mathcal{D}(\lambda)$ is concave and yields a lower
bound on the optimal value of \eqref{eq:objective fun} \cite{Boyd2010}.
Hence, the dual function $D(\lambda)$ is maximized with respect to
the variable $\lambda\in\mathbb{\mathbb{R}}^{m_{i}\times n}$, 
\begin{equation}
\mathcal{D}^{*}=\max_{\lambda\in\mathcal{Z}^{n}}\mathcal{D}(\lambda)=\max_{\lambda\in\mathcal{Z}^{n}}\mathcal{L}_{T}(x^{*},y^{*},\lambda);\label{eq:dual problem}
\end{equation}
we note that $\mathcal{Z}^{n}\subseteq\mathbb{\mathbb{R}}^{m_{i}\times n}$.
Slater condition guarantees zero duality gap and the existence of
a dual optimal solution $\lambda^{*}\in\mathcal{Z}^{n}$. When $(x^{*},y_{1}^{*},\cdots,y_{n}^{*})\in\chi\times Y^{n}$
solves the primal problem \eqref{eq:objective fun}-\eqref{eq:local linear constraints},
the primal and dual optimal vectors form a saddle-point for the Lagrangian
$\mathcal{L}_{T}$ \cite{Bertsekas1999}. Thus, based on the saddle
point theorem, if $w^{*}=(x^{*},y_{1}^{*},\cdots,y_{n}^{*},\lambda_{1}^{*},\cdots,\lambda_{n}^{*})\in\Omega$
is a saddle point for $\mathcal{L}_{T}$, then for all $w=(x,y_{1},\cdots,y_{n},\lambda_{1},\cdots,\lambda_{n})\in\chi\times Y^{n}\times\mathbb{\mathbb{R}}^{m_{i}\times n}=\Omega$,
we have \vspace{-0.2cm}
 
\begin{equation}
\mathcal{L}_{T}(x^{*},y^{*},\lambda)\leq\mathcal{L}_{T}(x^{*},y^{*},\lambda^{*})\leq\mathcal{L}_{T}(x,y,\lambda^{*}).\label{eq:FONC}
\end{equation}
Moreover, the Slater condition implies that the dual optimal set is
bounded; hence $\left\Vert \lambda_{i}^{*}\right\Vert \leq D_{\lambda}$
for all $i\in[n]$ for some finite $D_{\lambda}$ (see Lemma 3 in
\cite{Nedich2009}). A consequence of inequality \eqref{eq:FONC}
is that $\widetilde{w}=\left(\widetilde{x},\widetilde{y},\widetilde{\lambda}\right)\in\Omega$
approximately solves the primal problem with accuracy $\epsilon_{T}^{\mathcal{P}}\geq0$
if it satisfies, 
\[
0\leq\mathcal{L}_{T}(\widetilde{x},\widetilde{y},\lambda^{*})-\mathcal{L}_{T}(x^{*},y^{*},\lambda^{*})\leq\epsilon_{T}^{\mathcal{P}},
\]
that is, 
\begin{align}
0\leq\mathcal{L}_{T}(\widetilde{x},\widetilde{y},\lambda^{*})-\mathcal{P}^{*}\leq\epsilon_{T}^{\mathcal{P}}.\label{eq:lagrangian accuracy-primal}
\end{align}
Based on \eqref{eq:dual func}, the inequality \eqref{eq:lagrangian accuracy-primal}
can also be referred as dual feasibility. In addition, $\widetilde{w}=\left(\widetilde{x},\widetilde{y},\widetilde{\lambda}\right)\in\Omega$
approximately solves the dual problem with accuracy $\epsilon_{T}^{\mathcal{D}}\geq0$
if 
\[
0\leq\mathcal{L}_{T}(x^{*},y^{*},\lambda^{*})-\mathcal{L}_{T}(x^{*},y^{*},\widetilde{\lambda})\leq\epsilon_{T}^{\mathcal{D}},
\]
that is, 
\begin{align}
 & 0\leq\mathcal{D}^{*}-\mathcal{D}(\widetilde{\lambda})\leq\epsilon_{T}^{\mathcal{D}},\label{eq:lagrangian accuracy-dual}
\end{align}
which represents the dual sub-optimality. The conditions in \eqref{eq:lagrangian accuracy-primal}-\eqref{eq:lagrangian accuracy-dual}
can be combined to represent the duality gap as 
\[
\sum_{t=1}^{T}f_{t}^{\Delta}(\widetilde{w},w^{*})+\frac{1}{n}(\sum_{i=1}^{n}\phi_{i}^{\Delta}(\widetilde{w},w^{*})+H_{i}^{\Delta}(\widetilde{w},w^{*}))\leq\epsilon_{T}.
\]
where 
\begin{align*}
f_{t}^{\Delta}(w,w^{*}) & =f_{t}(x)-f_{t}(x^{*})\\
\phi_{i}^{\Delta}(w,w^{*}) & =\phi_{i}(y_{i})-\phi_{i}(y_{i}^{*})\\
H_{i}^{\Delta}(w,w^{*}) & =h_{1i}^{\Delta}(w,w^{*})+h_{2i}^{\Delta}(w,w^{*})\\
h_{1i}^{\Delta}(w,w^{*}) & =\left\langle x-x^{*},A_{i}^{T}\lambda_{i}^{*}\right\rangle +\left\langle \lambda_{i}-\lambda_{i}^{*},-r_{i}(x^{*},y_{i}^{*})\right\rangle \\
h_{2i}^{\Delta}(w,w^{*}) & =\left\langle y_{i}-y_{i}^{*},B_{i}^{T}\lambda_{i}^{*}\right\rangle ,
\end{align*}
and $\epsilon_{T}=\epsilon_{T}^{\mathcal{P}}+\epsilon_{T}^{\mathcal{D}}\geq0$.
This gap is illustrated in Figure \ref{fig:duality gap}. 

Analogous to the regret definition for O-ADMM algorithm \cite{Suzuki2013},
we can consider a sequence of decisions $w_{t}$, where $w_{t}\in\Omega$
for each $t$, instead of a fixed decision $\widetilde{w}$. Consequently,
the sequence $w_{t}$ approximately solves \eqref{eq:objective fun}
and \eqref{eq:local linear constraints} with accuracy $\epsilon_{T}$
if 
\begin{equation}
\sum_{t=1}^{T}f_{t}^{\Delta}(w_{t},w^{*})+\frac{1}{n}\left\{ \sum_{i=1}^{n}\phi_{i}^{\Delta}(w_{t},w^{*})+H_{i}^{\Delta}(w_{t},w^{*})\right\} \leq\epsilon_{T},\label{eq:epsilon accuracy}
\end{equation}
for the optimal solution $w^{*}\in\Omega$, referred to as \textit{fixed
case solutions} to distinguish them from the time-varying online solution
sequence $w_{t}$. Moreover, the mapping $H_{i}^{\Delta}(w,w^{*})$
can be expressed as 
\[
H_{i}^{\Delta}(w,w^{*})=\left\langle w_{i}(x)-w_{i}^{*}(x),H_{i}(w^{*})\right\rangle ,
\]
where $w_{i}(x)=\left[\begin{array}{ccc}
x & y_{i} & \lambda_{i}\end{array}\right]^{T}$, $w_{i}^{*}(x)=\left[\begin{array}{ccc}
x & y_{i}^{*} & \lambda_{i}^{*}\end{array}\right]^{T}$, and 
\[
H_{i}(w)=\left[\begin{array}{ccc}
0 & 0 & A_{i}^{T}\\
0 & 0 & B_{i}^{T}\\
-A_{i} & -B_{i} & 0
\end{array}\right]w_{i}(x)+\left[\begin{array}{c}
0\\
0\\
c_{i}
\end{array}\right].
\]
Since, the mapping $H_{i}(w)$ is affine in $w_{i}(x)$ and is defined
through a skew symmetric matrix, it is monotone, and consequently
\cite{Facchinei2003} 
\begin{align}
 & \left\langle w_{i}(x)-w_{i}^{*}(x),H_{i}(w)-H_{i}(w^{*})\right\rangle \geq0\nonumber \\
 & \left\langle w_{i}(x)-w_{i}^{*}(x),H_{i}(w)\right\rangle \geq\left\langle w_{i}(x)-w_{i}^{*}(x),H_{i}(w^{*})\right\rangle .\label{eq:monotone mapping}
\end{align}
Therefore, the inequality 
\begin{align}
\sum_{t=1}^{T}f_{t}^{\Delta}(w_{t},w^{*})+\frac{1}{n}\sum_{i=1}^{n}\left(\phi_{i}^{\Delta}(w_{t},w^{*})+\left\langle w_{i,t}(x)-w_{i}^{*}(x),H_{i}(w_{t})\right\rangle \right) & \leq\epsilon_{T}\label{eq:epsilonAccuracy2}
\end{align}
is a sufficient condition for \eqref{eq:epsilon accuracy}.

Finally, motivated by the inclusion of regularization terms in the
augmented Lagrangian method \cite{Bertsekas1999}, the term on the
left hand side of \eqref{eq:epsilonAccuracy2} is supplemented with
terms of the form $\frac{\rho}{2}||r_{i}(x_{i,t},y_{i,t})||^{2}$,
where $\rho>0$, to promote agents satisfy the local primal feasibility
constraints. In our setting, the sequence $w_{i,t}$ is constructed
from the distributed algorithm adopted by each agent $i$, specifically
$w_{t}(x_{i,t})=(x_{j,t},y_{t},\lambda_{t+1})\in\Omega$ at time $t$,
where $y_{t}=\left(y_{1,t},\cdots,y_{n,t}\right)$ and $\lambda_{t+1}=\left(\lambda_{1,t+1},\cdots,\lambda_{n,t+1}\right)$.
The social regret is thus defined as,\footnote{Note that this form of regret penalizes the deviation of each agent's local copy of the global variable from the best fixed global decision in hindsight.} 
$$R_{T}=\max_{j \in [n]} R_{j,T},$$ 
where
\begin{align}
R_{j,T} & =\sum_{t=1}^{T}f_{t}^{\Delta}(w_{t}(x_{j,t}),w^{*})+\frac{1}{n}\sum_{i=1}^{n}\left\{ \phi_{i}^{\Delta}(w_{t}(x_{j,t}),w^{*})+\left\langle w_{i,t}(x_{j,t})-w_{i}^{*}(x_{j,t}),H_{i}(w_{t}(x_{j,t}))\right\rangle +\frac{\rho}{2}||r_{i}(x_{i,t},y_{i,t})||^{2}\right\} .\label{eq:global regret}
\end{align}
Based on \eqref{eq:epsilonAccuracy2} we say that the sequence $w_{i,t}$
approximately solves \eqref{eq:objective fun} and \eqref{eq:local linear constraints}
with accuracy $\epsilon_{T}$ if it satisfies $R_{T}\leq\epsilon_{T}$.
Therefore, if the social regret is sub-linear with time, the online
algorithm performs as well as the best fixed case decision provided
with the complete sequence of cost functions a priori. In addition,
the sub-linearity of the social regret ensures that the local linear
constraints will be satisfied asymptotically.

\section{Main Result\label{sec:Main-Result}}

The main contribution of this paper is extending O-ADMM \cite{Suzuki2013a}
via Nesterov's Dual Averaging (DA) algorithm \cite{Nesterov2007}
and distributed subgradients (descent) method discussed in \cite{Yan2010,Nedic2009,SundharRam2010,Lobel2011},
to provide a distributed decision-making process for the optimization
problem discussed in \S\ref{sec:Problem-Statement} with a sub-linear
social regret; we refer to this procedure as online distributed ADMM
(OD-ADMM). The main challenge for the seamless integration of ADMM
with dual averaging and distributed gradient descent for OD-ADMM is
deriving and utilizing bounds on the network effect and the sub-optimality
of the local decisions on the quality of the agent's decisions on
the social regret. This objective is achieved by building on the existing
results reported in \cite{Yan2010,Duchi2012,Hosseini2013} regarding
the network contribution in distributed optimization, as well as extensions
of results discussed in~\cite{Wang2012a,Suzuki2013a,Hosseini2014}.
basic idea behind our convergence analysis is as follows.%
\footnote{All referenced lemmas are discussed in the Appendix.%
} First, in Lemmas \ref{lem:DA network effect} and \ref{lem:GD network effect},
we derive the gap between the local decisions and the average decision
over the network. Then, we provide the sub-optimality gap in Lemmas
\ref{lem:DA sub-optimality effect} and \ref{lem:GD sub-optimality effect}.
Finally, building on previous results, bounds on the social regret
are presented in Theorems \ref{thm:Regret(xi)-1} and \ref{thm:Regret(xi)_GD}.

The proposed algorithm updates the vector $(x_{i},y_{i},z_{i},\lambda_{i})$
for each agent $i\in[n]$ by alternately minimizing the Lagrangian
and augmented Lagrangian. In addition, the Lagrangian is linearized
based on network-level update, leading to a subgradient descent method
followed by a projection step onto the constraint set $\chi$. Specifically,
in the DA method, we let 
\[
z_{t+1}=z_{t}+\tilde{g}_{t},
\]
where $\tilde{g}_{t}=\nabla\mathcal{L}_{t}(x_{t})$, followed by 
\begin{equation}
x_{t+1}=\prod_{\chi}^{\psi}\left(z_{t+1},\alpha_{t}\right);\label{eq:centralized_x_update_DA}
\end{equation}
in this case, the parameter $\alpha_{t}$ is a non-increasing sequence
of positive functions and $\prod_{\chi}^{\psi}(\cdot)$ is the projection
operator onto $\chi$ defined as 
\begin{equation}
\prod_{\chi}^{\psi}(z_{t+1},\alpha_{t})\equiv\arg\min_{x\in\chi}\left\{ \left\langle z_{t+1},x\right\rangle +\frac{1}{\alpha_{t}}\psi(x)\right\} ;\label{eq:projection fun}
\end{equation}
the proximal function $\psi(x):\chi\rightarrow\mathbb{R}$ is continuously
differentiable and strongly convex . The inclusion of the proximal
function in the DA method as a regularizer prevents oscillations in
the projection step. 
to be strongly convex with respect to $\Vert.\Vert$, $\psi\geq0$,
and $\psi(0)=0$.

On the other hand, in the subgradient descent (GD) method, the aforementioned
steps in DA are replaced by 
\[
h_{t+1}=x_{t}-\alpha_{t}\tilde{g}_{t},
\]
followed by 
\begin{equation}
x_{t+1}=\prod_{\chi}h_{t+1}\equiv\arg\min_{x\in\chi}||x-h_{t+1}||.\label{eq:centralized_x_update_GD}
\end{equation}

Finally, the proposed online algorithm minimizes the augmented Lagrangian
over $y$ as 
\[
y_{t+1}=\arg\min_{y\in Y}\left\{ \mathcal{L}_{t}(x_{t+1},y,\lambda_{t+1})+\frac{\rho}{2}||r(x_{t+1},y)||^{2}\right\} ,
\]
and update the dual variable $\lambda$ as%
\footnote{Note that the index for the dual variable is one time step ahead of
the primal variables.%
} 
\[
\lambda_{t+2}=\lambda_{t+1}+\rho(Ax_{t+1}+By_{t+1}-c).
\]

The distributed algorithm can be considered as an approximate ADMM
by an agent $i$ via a convex combination of information provided
by its neighbors $N(i)$. Specifically, the global update step \eqref{eq:centralized_x_update_DA}
and \eqref{eq:centralized_x_update_GD} can be reformulated with a
distributed method. The underlying communication network can be represented
compactly as a doubly stochastic matrix $P\in\mathbb{R}^{n\times n}$
which preserves the zero structure of the Laplacian matrix $L(\mathcal{G})$.
For agents to have access to information contained in the subgradients
$\tilde{g}_{i,t}=\nabla\mathcal{L}_{i,t}(x_{i,t})$ there must be
information flow amongst the agents; as such, in our subsequent analysis
it will be assumed that the graph $\mathcal{G}$ is strongly connected.
A method to construct a doubly stochastic matrix $P$ of the required
form from the Laplacian of the network is provided in Proposition
\ref{prop:Strong}.

The online distributed ADMM (OD-ADMM) is presented in Algorithm \ref{alg:OD-ADMM}.
\begin{algorithm}
\SetAlgoLined 
 Initialize $x_{i,1}=0$ and $y_{i,1}=0$ for all $i=1,...,n$\\
 \For {$t=1$ \KwTo $T$}{ Adversary reveals $f_{t}(t)=\{f_{i,t}(t);\;\mbox{for}\;\forall i=1,...,n\}$\\
 Compute subgradient $g_{i}(t)\in\partial f_{i,t}(x_{i,t})$\\
 \ForEach {Agent $i$}{ $\lambda_{i,t+1}=\lambda_{i,t}+\rho(A_{i}x_{i,t}+B_{i}y_{i,t}-c_{i})$\ \label{eq:lambda_i update}

$x_{i,t+1}=H_{\alpha_{t}}(\lambda_{i,t+1},g_{i,t})$\ \label{eq:x_i update}

$r_{i}(x_{i,t+1},y)=A_{i}x_{i,t+1}+B_{i}y-c_{i}$\ \label{eq:r_i update}

$y_{i,t+1}=\mbox{argmin}_{y\in Y}\left\{ \phi_{i}(y)+\lambda_{i,t+1}^{T}r_{i}(x_{i,t+1},y)+\frac{\rho}{2}||r_{i}(x_{i,t+1},y)||^{2}\right\} $\ \label{eq:y_i update}

}} \protect\caption{Online Distributed ADMM (OD-ADMM)}

\label{alg:OD-ADMM} 
\end{algorithm}

The function $H_{\alpha_{t}}(\lambda_{i,t+1},g_{i,t})$ referred to
on line 8 of the algorithm represents a distributed update on the
primal variable. 
this paper, we consider two alternatives for this update. 

\subsection{OD-ADMM via Distributed Dual Averaging}

In this method, the dual sub-gradient at each node is updated as a
convex combination of its neighbor's dual subgradients and itself,
namely, 
\begin{equation}
z_{i,t+1}=\sum_{j=1}^{n}P_{ji}z_{j,t}+g_{i,t}+A_{i}^{T}\lambda_{i,t+1},\label{eq:DA local z update}
\end{equation}
and 
\begin{equation}
x_{i,t+1}=\prod_{\chi}^{\psi}(z_{i,t+1},\alpha_{t}),\label{eq:DA local x update}
\end{equation}
where the projection operator $\prod_{\chi}^{\psi}(\cdot)$ is defined
in \eqref{eq:projection fun}. Before presenting the convergence rate
of the proposed OD-ADMM algorithm we provide a few preliminary remarks
and definitions. Let us define the sequences of (network) average
dual subgradients $z_{t}$'s and average subgradients $g_{t}$'s as
\begin{align}
z_{t} & =\frac{1}{n}\sum_{i=1}^{n}z_{i,t},\;\;\; g_{t}=\frac{1}{n}\sum_{i=1}^{n}g_{i,t}.\label{eq:average z,g}
\end{align}
Thus, in the distributed DA method, the following update rule is introduced
similar to the standard DA algorithm, 
\begin{equation}
z_{t+1}=z_{t}+g_{t}+\frac{1}{n}\sum_{i=1}^{n}A_{i}^{T}\lambda_{i,t+1},\label{eq:z_bar update_DA}
\end{equation}
where the primal update is 
\begin{equation}
\theta_{t+1}=\Pi_{\chi}^{\psi}(z_{t+1},\alpha_{t}).\label{eq:theta update_DA}
\end{equation}

The regret analysis can now be presented as follows, where the intermediate
results required for its proof are relegated to the Appendix, namely
Lemmas \ref{lem:x_i-y-1}, \ref{lem:DA network effect}, and \ref{lem:DA sub-optimality effect}.
In particular, we show that with a proper choice of learning rate
in Lemmas \ref{lem:DA network effect} and \ref{lem:DA sub-optimality effect},
the network effect and the sub-optimality of average decisions are
sub-linear over time. Subsequently, building on these results, a sub-linear
regret bound for OD-ADMM using distributed DA method can be established
as formalized by the following result. 
\begin{thm}
\label{thm:Regret(xi)-1}Given the sequence $w_{i,t}$ generated by
Algorithm \ref{alg:OD-ADMM} where line \ref{eq:x_i update} applies
distributed dual-averaging method with $\psi(x^{*})\leq\Psi^{2}$
and $\alpha(t)=k/\sqrt{t}$, we have 
\begin{align}
R_{T}\leq J_{1}+J_{2}k\sqrt{T},\label{eq:regret_dist_bound-3}
\end{align}
where 
\[
J_{1}=\frac{D_{\lambda}}{\rho n}\sum_{i=1}^{n}\frac{\K_{i}}{\sigma_{1}(A_{i})},
\]
\begin{align*}
J_{2} & =2\mathcal{Q}(L_{f}+\Km)(\frac{2}{n}\sum_{i}(D_{\lambda}\sigma_{1}(A_{i})+2\zeta_{i}))
,
\end{align*}
with
\begin{align*}
 & \overline{\K}=\frac{1}{n}\sum_{i=1}^{n}\K_{i},\,\Km=\max_{i}\K_{i},\\
 & \K_{i}=\sqrt{m_{i}}L_{\phi}\frac{\sigma_{1}(A_{i})}{\sigma_{m_{i}}(B_{i}^{T})},\mbox{ and }\mathcal{Q}=\frac{\sqrt{n}}{1-\sigma_{2}(P)}.
\end{align*}
\end{thm}
\begin{IEEEproof}
Based on the definition of $f_{t}$ we have 
\[
f_{t}^{\Delta}(w_{t}(x_{j,t}),w^{*})=\frac{1}{n}\sum_{i=1}^{n}f_{i,t}^{\Delta}(w_{t}(x_{j,t}),w^{*}),
\]
where $w_{t}(x_{j,t})=(x_{j,t},y_{t},\lambda_{t+1})\in\Omega$ and
thus 
\[
R_{j,T}=\frac{1}{n}\sum_{i,t}f_{i,t}^{\Delta}(w_{t}(x_{j,t}),w^{*})+\phi_{i}^{\Delta}(w_{t}(x_{j,t}),w^{*})+H_{i}^{\Delta}(w_{t}(x_{j,t}),w^{*})+\frac{\rho}{2}||r_{i}(x_{i,t},y_{i,t})||^{2}.
\]

In the meantime as $f_{t}$ is $L$-Lipschitz and convex, we have
\begin{align}
f_{t}^{\Delta}(w_{t}(x_{j,t}),w^{*}) & =f_{t}(x_{j,t})-f_{t}(\theta_{t})+f_{t}(\theta_{t})-f_{t}(x^{*})\nonumber \\
 & \leq L_{f}\left\Vert x_{j,t}-\theta_{t}\right\Vert +\left\langle g_{t},\theta_{t}-x^{*}\right\rangle .\label{eq:ft(wt,w)-1}
\end{align}
The first term in \eqref{eq:ft(wt,w)-1} represents the network effect
in the regret bound, i.e., the deviation of local primal variable
at each node from the average primal variable. Lemma \ref{lem:DA network effect}
in the Appendix on the other hand, provides a bound on the network
effect using the DA method. Therefore, replacing line \ref{eq:x_i update}
of Algorithm \ref{alg:OD-ADMM} with the distributed DA method implies
that 
\begin{equation}
\Vert\theta_{t}-x_{j,t}\Vert\leq\alpha_{t-1}\frac{\sqrt{n}\,(L_{f}+\Km)}{1-\sigma_{2}(P)}.\label{eq:conv_5-2}
\end{equation}
Moreover, from the integral test with $\alpha_{t}=k/\sqrt{t}$ it
follows that%
\footnote{Note that $\frac{1}{\sqrt{t}}$ is a non increasing positive function
and the integral test leads to $\sum_{t=1}^{T}\frac{1}{\sqrt{t}}\leq2\sqrt{T}-1$.%
} 
\begin{equation}
\sum_{t=1}^{T}\alpha_{t-1}\leq2k\sqrt{T}.\label{eq:sum_alpha-1}
\end{equation}
Hence, from \eqref{eq:conv_5-2} and \eqref{eq:sum_alpha-1} it follows
that 
\begin{equation}
\sum_{t=1}^{T}\Vert x_{j,t}-\theta_{t}\Vert\leq2k\sqrt{T}\mathcal{Q}(L_{f}+\Km).\label{eq:conv_6-1}
\end{equation}

The second term in \eqref{eq:ft(wt,w)-1} represents the sub-optimality
of the procedure due to using the subgradient method. Applying Lemma
\ref{lem:DA sub-optimality effect} (Appendix) with \eqref{eq:z_bar update_DA}
and \eqref{eq:theta update_DA} implies that 
\begin{align}
 & \sum_{t=1}^{T}\left\langle g_{t},\theta_{t}-x^{*}\right\rangle \leq\sum_{t=1}^{T}[\frac{\alpha_{t}}{2}\left\Vert \vphantom{g_{t+1}}\right.g_{t+1}+\frac{1}{n}\sum_{i=1}^{n}A_{i}^{T}\lambda_{i,t+2}\left.\vphantom{g_{t+1}}\right\Vert _{*}^{2}\nonumber \\
 & \hphantom{\leq}\;\;\;\;\;\;\;\;\;\;\;\;\;\;\;\;\;\;\;\;\;\;+\frac{1}{n}\sum_{i=1}^{n}\left\langle \lambda_{i,t+1},A_{i}(x^{*}-\theta_{t})\right\rangle ]+\frac{1}{\alpha_{T}}\psi\left(x^{*}\right).\label{eq:f_3-1}
\end{align}

The first term on the the right hand side of \eqref{eq:f_3-1} represents
the sub-optimality of Lagrangian function, defined in \eqref{eq: Lagrangian},
with respect to the global variable $x$ and is bounded as%
\footnote{Note that $\Vert Qx\Vert\leq\sigma_{1}(Q)||x||$ for any matrix $Q\in\mathbb{R}^{m\times n}$
and vector $x\in\mathbb{R}^{n}$.%
} 
\begin{align}
\sum_{t=1}^{T}\frac{\alpha_{t}}{2}\left\Vert \vphantom{g_{t+1}}\right.g_{t+1}+\frac{1}{n}\sum_{i=1}^{n}A_{i}^{T}\lambda_{i,t+2}\left.\vphantom{g_{t+1}}\right\Vert _{*}^{2}\leq%
(\max_{t}\left\Vert g_{t+1}\right\Vert _{*}+\frac{1}{n}\sum_{i=1}^{n}\sigma_{1}(A_{i})\max_{t}\left\Vert \lambda_{i,t}\right\Vert _{*})^{2}\sum_{t=1}^{T}\frac{\alpha_{t}}{2}.\label{eq:norm(g+A*lambda)-1}
\end{align}

We now proceed to bound the individual terms in \eqref{eq:norm(g+A*lambda)-1}.
By optimality of line \ref{eq:y_i update} in Algorithm \ref{alg:OD-ADMM}
and applying line \ref{eq:lambda_i update}, we have 
\[
\nabla_{y}\phi_{i}(y_{i,t})=-B_{i}^{T}\left(\lambda_{i,t}+\rho r_{i}(x_{i,t},y_{i,t})\right)=-B_{i}^{T}\lambda_{i,t+1},
\]
for all $i\in[n]$ and $t\in[T]$. Moreover, since $||\nabla_{y}\phi_{i}(y_{i,t})||\leq L_{\phi}$,
we have $||B_{i}^{T}\lambda_{i,t+1}||\leq L_{\phi}$. Thus, $\lambda_{i,t}$
is bounded as 
\begin{align}
||\lambda_{i,t}|| & \leq||(B_{i}B_{i}^{T})^{-1}B_{i}||_{F\;}||B_{i}^{T}\lambda_{i}||\leq L_{\phi}(\sum_{j=1}^{m_{i}}\frac{1}{\sigma_{j}^{2}(B_{i}^{T})})^{1/2}\leq\frac{\sqrt{m_{i}}L_{\phi}}{\sigma_{m_{i}}(B_{i}^{T})},\label{eq:norm(lambda)-1}
\end{align}
which implies that $||A_{i}^{T}\lambda_{i,t}||\leq\sqrt{m_{i}}L_{\phi}\sigma_{1}(A_{i})/\sigma_{m_{i}}(B_{i}^{T}).$
Based on Lipschitz continuity of $f_{t}$, we have that $\left\Vert g_{t+1}\right\Vert _{*}\leq L_{f},$
and subsequently \eqref{eq:norm(g+A*lambda)-1} is bounded as 
\begin{equation}
\sum_{t=1}^{T}\frac{\alpha_{t}}{2}\left\Vert \vphantom{g_{t+1}}\right.g_{t+1}+\frac{1}{n}\sum_{i=1}^{n}A_{i}^{T}\lambda_{i,t+2}\left.\vphantom{g_{t+1}}\right\Vert _{*}^{2}\leq(L_{f}+\overline{\K})^{2}k\sqrt{T}.\label{eq:norm(g+A*lambda)) bound-1}
\end{equation}
The second term in the inequality \eqref{eq:f_3-1} represents the
sub-optimality of centralized decision $\theta$ with respect to the
linear constraints. In order to analyze this term, it is first expanded
into two terms representing the sub-optimality of local decision $x_{i,t}$
and the network effect, respectively, i.e., 
\begin{align}
 & \left\langle \lambda_{i,t+1},A_{i}(x^{*}-\theta_{t})\right\rangle \nonumber \\
 & =\left\langle \lambda_{i,t+1},A_{i}(x^{*}-x_{j,t})\right\rangle +\left\langle \lambda_{i,t+1},A_{i}(x_{j,t}-\theta_{t})\right\rangle \nonumber \\
 & =\left\langle \lambda_{i,t+1},A_{i}(x^{*}-x_{j,t})\right\rangle +\left\langle \lambda_{i}^{*}-\lambda_{i,t+1},-r_{i}(x_{j,t},y_{i,t})\right\rangle \nonumber \\
 & \;\;\;\;+\left\langle \lambda_{i,t+1},A_{i}(x_{j,t}-\theta_{t})\right\rangle +\left\langle \lambda_{i}^{*}-\lambda_{i,t+1},r_{i}(x_{j,t},y_{i,t})\right\rangle \nonumber \\
 & =-h_{1i}^{\Delta}(w_{t}(x_{j,t}),w^{*})+\left\langle \lambda_{i,t+1},A_{i}(x_{j,t}-\theta_{t})\right\rangle +\left\langle \lambda_{i}^{*}-\lambda_{i,t+1},r_{i}(x_{j,t},y_{i,t})\right\rangle .\label{eq:lambdaA(x-theta)-1}
\end{align}
Based on the network effect introduced in \eqref{eq:conv_5-2}, we
have 
\begin{equation}
\left\langle \lambda_{i,t+1},A_{i}(x_{j,t}-\theta_{t})\right\rangle \leq\sigma_{1}(A_{i})||\lambda_{i,t+1}||\;\Vert x_{j,t}-\theta_{t}\Vert.\label{eq:<lambda,A(x-theta)> bound-1}
\end{equation}
Moreover, applying \eqref{eq:conv_6-1} and \eqref{eq:norm(lambda)-1}
to \eqref{eq:<lambda,A(x-theta)> bound-1}, it follows that 
\begin{align}
 & \sum_{t=1}^{T}\left\langle \lambda_{i,t+1},A_{i}(x_{j,t}-\theta_{t})\right\rangle \leq2k\sqrt{T}\K_{i}\mathcal{Q}(L_{f}+\Km).\label{eq:sum(lambda,A(x-theta))-1}
\end{align}
The final term in \eqref{eq:lambdaA(x-theta)-1} represents the first
order necessary condition for optimality of the dual problem at $\lambda_{i,t+1}$.
By applying line 7 of the algorithm and an inner product equality,
we obtain%
\footnote{Namely, using the identity $\left\langle v_{1}-v_{2},v_{3}+v_{4}\right\rangle =\frac{1}{2}(||v_{4}-v_{2}||^{2}-||v_{4}-v_{1}||^{2}+||v_{3}+v_{1}||^{2}-||v_{3}+v_{2}||^{2}$.%
} 
\begin{align}
\left\langle \lambda_{i}^{*}-\lambda_{i,t+1},r_{i}(x_{j,t},y_{i,t})\right\rangle  & =\frac{1}{\rho}\left\langle \lambda_{i}^{*}-\lambda_{i,t+1},\lambda_{i,t+1}-\lambda_{i,t}\right\rangle +\left\langle \lambda_{i}^{*}-\lambda_{i,t+1},A_{i}(x_{j,t}-x_{i,t})\right\rangle \nonumber \\
 & =\frac{1}{2\rho}(-\left\Vert \lambda_{i,t+1}-\lambda_{i,t}\right\Vert ^{2}+\left\Vert \lambda_{i}^{*}-\lambda_{i,t}\right\Vert ^{2}-\left\Vert \lambda_{i}^{*}-\lambda_{i,t+1}\right\Vert ^{2})\nonumber \\
 & \;\;\;\;+\left\langle \lambda_{i}^{*}-\lambda_{i,t+1},A_{i}(x_{j,t}-x_{i,t})\right\rangle \nonumber \\
 & =\frac{1}{2\rho}(\left\Vert \lambda_{i}^{*}-\lambda_{i,t}\right\Vert ^{2}-\left\Vert \lambda_{i}^{*}-\lambda_{i,t+1}\right\Vert ^{2})-\frac{\rho}{2}\left\Vert r_{i}(x_{i,t},y_{i,t})\right\Vert ^{2}\nonumber \\
 & \;\;\;\;+\left\langle \lambda_{i}^{*}-\lambda_{i,t+1},A_{i}(x_{j,t}-x_{i,t})\right\rangle .\label{eq:<lambda*-lambda,ri>}
\end{align}
Resolving the telescoping sum 
\[
\sum_{t=1}^{T}\left\Vert \lambda_{i}^{*}-\lambda_{i,t}\right\Vert ^{2}-\left\Vert \lambda_{i}^{*}-\lambda_{i,t+1}\right\Vert ^{2},
\]
using the fact $\lambda_{i,1}=0$, it now follows that 
\begin{align*}
\sum_{t=1}^{T}\left\langle \lambda_{i}^{*}-\lambda_{i,t+1},r_{i}(x_{i,t},y_{i,t})\right\rangle  & \leq\frac{1}{2\rho}(\left\Vert \lambda_{i}^{*}\right\Vert ^{2}-\left\Vert \lambda_{i}^{*}-\lambda_{i,T+1}\right\Vert ^{2})-\frac{\rho}{2}\sum_{t=1}^{T}\left\Vert r_{i}(x_{i,t},y_{i,t})\right\Vert ^{2}\\
 & \leq\frac{1}{2\rho}(2\left\Vert \lambda_{i}^{*}\right\Vert \left\Vert \lambda_{i,T+1}\right\Vert )-\frac{\rho}{2}\sum_{t=1}^{T}\left\Vert r_{i}(x_{i,t},y_{i,t})\right\Vert ^{2}.
\end{align*}
Applying \eqref{eq:norm(lambda)-1} in conjunction with the assumption
$\left\Vert \lambda_{i}^{*}\right\Vert \leq D_{\lambda}$, 
\begin{align}
 & \sum_{t=1}^{T}\left\langle \lambda_{i}^{*}-\lambda_{i,t+1},r_{i}(x_{i,t},y_{i,t})\right\rangle \leq\frac{D_{\lambda}\K_{i}}{\rho\sigma_{1}(A_{i})}-\frac{\rho}{2}\sum_{t=1}^{T}\left\Vert r_{i}(x_{i,t},y_{i,t})\right\Vert ^{2}.\label{eq:sum(lambda-lambda,r)-1}
\end{align}
The last term in \eqref{eq:<lambda*-lambda,ri>} can also be bounded
as 
\begin{align}
\left\langle \lambda_{i}^{*}-\lambda_{i,t+1},A_{i}(x_{j,t}-x_{i,t})\right\rangle  & \leq(D_{\lambda}\sigma_{1}(A_{i})+\zeta_{i})\Vert x_{j,t}-\theta_{t}\Vert\Vert x_{i,t}-\theta_{t}\Vert\nonumber \\
 & \leq2\alpha_{t-1}\mathcal{Q}(L_{f}+\Km)(D_{\lambda}\sigma_{1}(A_{i})+\zeta_{i})\label{eq: <lambda*-lambda,Ai(xj-xi)>}
\end{align}

Substituting \eqref{eq:sum(lambda,A(x-theta))-1}, \eqref{eq:sum(lambda-lambda,r)-1},
and \eqref{eq: <lambda*-lambda,Ai(xj-xi)>} into \eqref{eq:lambdaA(x-theta)-1},
\begin{align}
\frac{1}{n}\sum_{i=1}^{n}\sum_{t=1}^{T}\left\langle \lambda_{i,t+1},A_{i}\left(x^{*}-\theta_{t}\right)\right\rangle \leq J_{1}-\frac{1}{n}\sum_{i,t}\left\{ h_{1i}^{\Delta}(w_{t}(x_{j,t}),w^{*})+\frac{\rho}{2}||r_{i}(x_{i,t},y_{i,t})||^{2}\right\} \nonumber \\
+4k\sqrt{T}\mathcal{Q}(L_{f}+\Km)(\frac{1}{n}\sum_{i=1}^{n}D_{\lambda}\sigma_{1}(A_{i})+\overline{\mathcal{\K}}).\label{eq:sum(lambda,A(x-theta))bound-1}
\end{align}
Applying $\psi\left(x^{*}\right)\leq\Psi^{2}$, $\alpha_{T}=k/\sqrt{T}$,
and substituting \eqref{eq:norm(g+A*lambda)) bound-1}, \eqref{eq:sum(lambda,A(x-theta))bound-1}
into \eqref{eq:f_3-1} and simplifying, the sub-optimality of the primal
problem at the global decision $\theta_{t}$ can now be represented
as 
\begin{align}
\sum_{t=1}^{T}\left\langle g_{t},\theta_{t}-x^{*}\right\rangle \leq J_{1}-\frac{1}{n}\sum_{i,t}\left\{ h_{1i}^{\Delta}(w_{t}(x_{j,t}),w^{*})+\frac{\rho}{2}||r_{i}(x_{i,t},y_{i,t})||^{2}\right\} \nonumber \\
+k\sqrt{T}\times(2\mathcal{Q}(L_{f}+\Km)(\frac{2}{n}\sum_{i}(D_{\lambda}\sigma_{1}(A_{i})+2\zeta_{i})).
\label{eq:sum(g,theta-x)-1}
\end{align}
Based on our assumption of convexity of $\phi_{i}(\cdot)$, we have
\begin{align}
\phi_{i}^{\Delta}(w_{i,t},w^{*}) & \leq\left\langle \nabla_{y}\phi_{i}(y_{i,t}),y_{i,t}-y_{i}^{*}\right\rangle \nonumber \\
 & \leq-\left\langle B_{i}^{T}\lambda_{i,t+1},y_{i,t}-y_{i}^{*}\right\rangle \nonumber \\
 & =-h_{2i}^{\Delta}(w_{t}(x_{j,t}),w^{*}).\label{eq:phi(wt,w)-1}
\end{align}
Combining \eqref{eq:sum(g,theta-x)-1} and \eqref{eq:conv_6-1} into
\eqref{eq:ft(wt,w)-1} and adding \eqref{eq:phi(wt,w)-1}, the regret
can thereby be bounded as 
\[
R_{j,T}\leq J_{1}+J_{2}k\sqrt{T}
\]
for all $j \in [n]$ and thus the social regret is bounded as $R_{T}\leq J_{1}+J_{2}k\sqrt{T}$. Note that the social regret represents the worst case 
regret amongst the agents in the network. 
\end{IEEEproof}
The above theorem validates the ``good'' performance of OD-ADMM
via dual averaging by demonstrating a sub-linear social regret. In
addition, this social regret highlights the importance of the underlying
interaction topology through $\sigma_{2}(P)$ and certain condition
measure of local linear constraints through $\sigma_{1}(A_{i})$ and
$\sigma_{m_{i}}(B_{i})$. A well known measure of network connectivity
is the second smallest eigenvalue of the graph Laplacian $L(\mathcal{G})$
denoted by $\Lambda_{2}(\mathcal{G})$. Since the communication matrix
$P$ is formed as proposed in Proposition \ref{prop:Strong}, $1-\sigma_{2}(P)$
is proportional to $\Lambda_{2}(\mathcal{G})$ implying that high
network connectivity promotes good performance of the proposed OD-ADMM
algorithm with the embedded distributed implementation of dual averaging.

\subsection{OD-ADMM via Distributed Gradient Descent}

In this section, the local primal variable is updated using distributed
GD method. In this method $x_{i,t}$ is updated as a convex combination
of its neighbor's local primal variables and itself, moving in the
direction of decreasing the Lagrangian function, 
\begin{equation}
h_{i,t+1}=\sum_{j=1}^{n}P_{ji}x_{j,t}-\alpha_{t}(g_{i,t}+A_{i}^{T}\lambda_{i,t+1}),\label{eq:GD local h update}
\end{equation}
followed by the projection onto the convex set $\chi$, 
\begin{equation}
x_{i,t+1}=\prod_{\chi}h_{i,t+1}.\label{eq:GD local x update}
\end{equation}

In the distributed GD method, we first define the (network) average
primal variable as 
\begin{equation}
\theta_{t}=\frac{1}{n}\sum_{i=1}^{n}x_{i,t}.\label{eq:theta update_GD}
\end{equation}

The regret analysis for the OD-ADMM via the GD method can now be presented
as follows, where the intermediate results required for its proof
are relegated to the Appendix, namely Proposition \ref{prop:orthogonal proj},
and Lemmas \ref{lem:GD network effect} and \ref{lem:GD sub-optimality effect}.
In particular, with a proper choice of the learning rate in Lemmas
\ref{lem:GD network effect} and \ref{lem:GD sub-optimality effect},
we can show that the network effect and sub-optimality of the average
decision are sub-linear over time. Then, building on these results,
a sub-linear social regret bound for OD-ADMM using distributed GD
method can be obtained; this is formalized in the following theorem. 
\begin{thm}
\label{thm:Regret(xi)_GD}Given the sequence $w_{i,t}$ generated
by Algorithm \ref{alg:OD-ADMM}, where line \ref{eq:x_i update} applies
the distributed GD method with $\alpha(t)=k/\sqrt{t}$, we have 
\begin{align}
R_{T}\leq J_{1}+J_{2}k\sqrt{T},\label{eq:regret_dist_bound}
\end{align}
where 
\[
J_{1}=\frac{D_{\lambda}}{\rho n}\sum_{i=1}^{n}\frac{\K_{i}}{\sigma_{1}(A_{i})}+\frac{D_{\chi}^{2}}{2k},
\]
\begin{align*}
J_{2} & =4\mathcal{Q}(L_{f}+\Km)(\frac{1}{n}\sum_{i=1}^{n}D_{\lambda}\sigma_{1}(A_{i})+2\overline{\K})+2(L_{f}+\overline{\K})^{2}+8L_{f}\mathcal{Q}(L_{f}+\overline{\K}),
\end{align*}
with 
\begin{align*}
 & \overline{\K}=\frac{1}{n}\sum_{i=1}^{n}\K_{i},\,\Km=\max_{i}\K_{i},\\
 & \K_{i}=\sqrt{m_{i}}L_{\phi}\frac{\sigma_{1}(A_{i})}{\sigma_{m_{i}}(B_{i}^{T})},\mbox{ and }\mathcal{Q}=\frac{\sqrt{n}}{1-\sigma_{2}(P)}.
\end{align*}
\end{thm}
\begin{IEEEproof}
Based on the definition of $f_{t}$ we have 
\[
f_{t}^{\Delta}(w_{t}(x_{j,t}),w^{*})=\frac{1}{n}\sum_{i=1}^{n}f_{i,t}^{\Delta}(w_{t}(x_{j,t}),w^{*}),
\]
where $w_{t}(x_{j,t})=(x_{j,t},y_{t},\lambda_{t+1})\in\Omega$ and
thus 
\[
R_{j,T}=\frac{1}{n}\sum_{i,t}f_{i,t}^{\Delta}(w_{t}(x_{j,t}),w^{*})+\phi_{i}^{\Delta}(w_{t}(x_{j,t}),w^{*})+H_{i}^{\Delta}(w_{t}(x_{j,t}),w^{*})+\frac{\rho}{2}||r_{i}(x_{i,t},y_{i,t})||^{2}.
\]

As $f_{t}$ is $L$-Lipschitz and convex, we have 
\begin{align}
f_{t}^{\Delta}(w_{t}(x_{j,t}),w^{*}) & =f_{t}(x_{j,t})-f_{t}(\theta_{t})+f_{t}(\theta_{t})-f_{t}(x^{*})\nonumber \\
 & \leq L_{f}\left\Vert x_{j,t}-\theta_{t}\right\Vert +\left\langle g_{t},\theta_{t}-x^{*}\right\rangle .\label{eq:ft(wt,w)}
\end{align}
The first term in \eqref{eq:ft(wt,w)} represents the network effect
in the regret bound, i.e., the deviation of local primal variable
at each node from the average primal variable. In the meantime, Lemmas
\ref{lem:DA network effect} and \ref{lem:GD network effect} in the
Appendix provide bounds on the network effect when the DA and GD methods
are used in line \ref{eq:x_i update} of Algorithm \ref{alg:OD-ADMM},
respectively. Therefore, replacing line \ref{eq:x_i update} of Algorithm
\ref{alg:OD-ADMM} with the distributed DA method implies 
\begin{equation}
\Vert\theta_{t}-x_{j,t}\Vert\leq\alpha_{t-1}\frac{\sqrt{n}(L_{f}+\Km)}{1-\sigma_{2}(P)},\label{eq:conv_5}
\end{equation}
and with the distributed GD, 
\begin{equation}
\Vert\theta_{t}-x_{j,t}\Vert\leq2\sqrt{n}(L_{f}+\Km)\sum_{k=1}^{t-1}\alpha_{t-k}\sigma_{2}(P)^{k-1}.\label{eq:conv_5-1}
\end{equation}
Moreover, from \eqref{eq:conv_5} (and \eqref{eq:sum_alpha-1}) it
follows that 
\begin{equation}
\sum_{t=1}^{T}\Vert x_{j,t}-\theta_{t}\Vert\leq4k\sqrt{T}(\mathcal{Q}(L_{f}+\Km).\label{eq:conv_6}
\end{equation}
Note that the upper bound in \eqref{eq:conv_6} is more conservative
than the distributed DA method by a factor of $2$.

The second term in \eqref{eq:ft(wt,w)} represents the sub-optimality
due to using sub-gradient method. Applying Lemma \ref{lem:GD sub-optimality effect}
with \eqref{eq:GD local h update}- \eqref{eq:theta update_GD} then
leads to 
\begin{align}
\sum_{t=1}^{T}\left\langle g_{t},\theta_{t}-x^{*}\right\rangle  & \leq\frac{2}{n^{2}}\sum_{t=1}^{T}\alpha_{t}(\sum_{i=1}^{n}||g_{i,t}+A_{j}^{T}\lambda_{i,t+1}||)^{2}+\sum_{t=1}^{T}(4L_{f}\sum_{k=0}^{t-1}\alpha_{t-k}\sigma_{2}(P)^{k})\frac{1}{n}\sum_{i=1}^{n}||(g_{i,t}+A_{i}^{T}\lambda_{i,t+1}||\nonumber \\
 & +\frac{1}{n}\sum_{t=1}^{T}\sum_{i=1}^{n}\left\langle A_{i}^{T}\lambda_{i,t+1},x^{*}-\theta_{t}\right\rangle +\frac{1}{2\alpha_{1}}D_{\chi}^{2}.\label{DUPLICATE: eq:f_3}
\end{align}

The first term on the right hand side of \eqref{DUPLICATE: eq:f_3}
is bounded as 
\begin{align}
 & \frac{2}{n^{2}}\sum_{t=1}^{T}\alpha_{t}(\sum_{i=1}^{n}||g_{i,t}+A_{i}^{T}\lambda_{i,t+1}||)^{2}\nonumber \\
 & \leq2(\frac{1}{n}\sum_{i=1}^{n}\max_{t}\left\Vert g_{t+1}\right\Vert +\frac{1}{n}\sum_{i=1}^{n}\sigma_{1}(A_{i})\max_{t}\left\Vert \lambda_{i,t}\right\Vert )^{2}\sum_{t=1}^{T}\sum_{t=1}^{T}\alpha_{t}\nonumber \\
 & \leq2(L_{f}+\overline{\K})^{2}k\sqrt{T}.\label{eq:norm(g+A*lambda)}
\end{align}

Similarly, the second term on the right hand side of \eqref{DUPLICATE: eq:f_3}
is bounded as 
\begin{align}
 & \sum_{t=1}^{T}(4L_{f}\sum_{k=0}^{t-1}\alpha_{t-k}\sigma_{2}(P)^{k})\frac{1}{n}\sum_{i=1}^{n}||(g_{i,t}+A_{i}^{T}\lambda_{i,t+1}||\nonumber \\
 & \leq8L_{f}k\sqrt{T}\mathcal{Q}(L_{f}+\overline{\K}).\label{eq:norm(g+A*lambda)2}
\end{align}
We now proceed to bound the third term in \eqref{DUPLICATE: eq:f_3}
and from \eqref{eq:lambdaA(x-theta)-1} we have 
\begin{align}
\left\langle \lambda_{i,t+1},A_{i}(x^{*}-\theta_{t})\right\rangle  & =-h_{1i}^{\Delta}(w_{t}(x_{j,t}),w^{*})+\left\langle \lambda_{i,t+1},A_{i}(x_{j,t}-\theta_{t})\right\rangle +\left\langle \lambda_{i}^{*}-\lambda_{i,t+1},r_{i}(x_{j,t},y_{i,t})\right\rangle .\label{eq:lambdaA(x-theta)}
\end{align}
Analogous to the proof of Theorem \ref{thm:Regret(xi)-1}, the second
term of \eqref{eq:lambdaA(x-theta)} is bounded as 
\begin{align}
 & \sum_{t=1}^{T}\left\langle \lambda_{i,t+1},A_{i}(x_{j,t}-\theta_{t})\right\rangle \leq4k\sqrt{T}\K_{i}\mathcal{Q}(L_{f}+\Km).\label{eq:sum(lambda,A(x-theta))}
\end{align}
Based on \eqref{eq:sum(lambda-lambda,r)-1} and \eqref{eq: <lambda*-lambda,Ai(xj-xi)>},
we have 
\begin{align}
 & \sum_{t=1}^{T}\left\langle \lambda_{i}^{*}-\lambda_{i,t+1},r_{i}(x_{j,t},y_{i,t})\right\rangle \leq\frac{D_{\lambda}\K_{i}}{\rho\sigma_{1}(A_{i})}-\frac{\rho}{2}\sum_{t=1}^{T}\left\Vert r_{i}(x_{i,t},y_{i,t})\right\Vert ^{2}+4k\sqrt{T}\mathcal{Q}(L_{f}+\Km)(D_{\lambda}\sigma_{1}(A_{i})+\zeta_{i}).\label{eq:sum(lambda-lambda,r)}
\end{align}
Substituting \eqref{eq:sum(lambda,A(x-theta))} and \eqref{eq:sum(lambda-lambda,r)}
into \eqref{eq:lambdaA(x-theta)}, 
\begin{align}
\frac{1}{n}\sum_{i=1}^{n}\sum_{t=1}^{T}\left\langle \lambda_{i,t+1},A_{i}\left(x^{*}-\theta_{t}\right)\right\rangle \leq & -\frac{1}{n}\sum_{i,t}[h_{1i}^{\Delta}(w_{t}(x_{j,t}),w^{*})+\frac{\rho}{2}||r_{i}(x_{i,t},y_{i,t})||^{2}]\nonumber \\
 & +4k\sqrt{T}\mathcal{Q}(L_{f}+\Km)(\frac{1}{n}\sum_{i=1}^{n}D_{\lambda}\sigma_{1}(A_{i})+2\overline{\K})+\frac{D_{\lambda}}{\rho n}\sum_{i=1}^{n}\frac{\K_{i}}{\sigma_{1}(A_{i})}.\label{eq:sum(lambda,A(x-theta))bound}
\end{align}
From \eqref{eq:norm(g+A*lambda)}, \eqref{eq:norm(g+A*lambda)2},
and \eqref{eq:sum(lambda,A(x-theta))bound}, the bound on \eqref{DUPLICATE: eq:f_3}
for $\alpha_{t}=k/\sqrt{t}$ is simplified to 
\begin{align}
\sum_{t=1}^{T}\sum_{t=1}^{T}\left\langle g_{t},\theta_{t}-x^{*}\right\rangle \leq & J_{1}-\frac{1}{n}\sum_{i,t}[h_{1i}^{\Delta}(w_{t}(x_{j,t}),w^{*})+\frac{\rho}{2}||r_{i}(x_{i,t},y_{i,t})||^{2}]\nonumber \\
 & +2k\sqrt{T}(\mathcal{Q}(L_{f}+\Km)(\frac{2}{n}\sum_{i=1}^{n}D_{\lambda}\sigma_{1}(A_{i})+4\overline{\K})+(L_{f}+\overline{\K})^{2}+4L_{f}\mathcal{Q}(L_{f}+\overline{\K})).\label{eq:sum(g,theta-x)}
\end{align}
In the meantime, based on the convexity of $\phi_{i}$, 
\begin{align}
\phi_{i}^{\Delta}(w_{t},w^{*}) & \leq-h_{2i}^{\Delta}(w_{t}(x_{j,t}),w^{*}).\label{eq:phi(wt,w)}
\end{align}
Combining \eqref{eq:sum(g,theta-x)} and \eqref{eq:conv_6} into \eqref{eq:ft(wt,w)}
and adding \eqref{eq:phi(wt,w)} we obtain, 
\[
R_{j,T}\leq J_{1}+J_{2}k\sqrt{T}.
\]
for all $j \in [n]$ and thus $R_{T}\leq J_{1}+J_{2}k\sqrt{T}$.
\end{IEEEproof}
Theorems \ref{thm:Regret(xi)-1} and \ref{thm:Regret(xi)_GD} and their
proofs provide a basis for comparing two effective methods for evaluating
OD-ADMM.
The bounds provided in \eqref{eq:conv_6-1} and \eqref{eq:conv_6} for the
distributed DA and GD, respectively, although conservative, hint at the fact that
in the distributed DA, the local copies of the global variable $x\in\chi$ might converge faster to consensus
in the worst case scenario. 
Moreover, as discussed in the introduction, the distributed DA approach does not require
sharing the primal variables $x_{i,t}$ amongst the agents, preserving their
privacy during the distributed decision making process; this feature of the DA approach however is not
shared by embedding the distributed GD in OD-ADMM.

\section{Example - Formation Acquisition with Points of Interest and Boundary
Constraints \label{sec:Example---Distributed}}
In order to demonstrate the applicability of the developed OD-ADMM, we consider
the following problem from the area of distributed robotics.
This so-called formation acquisition problem is 
as follows. Consider $n$ planar robots (agents) where
the position of agent $i$, denoted as $y_{i}$, is restricted to
the convex set $Y=\left[-1,1\right]^{2}$. 
The centroid of the formation
is denoted by $x\in\mathbb{R}^{2}$ and similarly constrained to $\mathcal{X}=Y$.
The formation shape is defined for each agent by its offset $c_{i}$
from the centroid, namely $x-y_{i}=c_{i}$. 
There is a known boundary
$S$ which agents are required to avoid by increasing the distance
to the boundary $\mbox{dist}(y_{i},S)=\inf_{x\in S}\left\Vert x-y_{i}\right\Vert $.
This is achieved with a penalty function $\phi_{i}(y_{i})=\left(\mbox{dist}(y_{i},S)+1\right)^{-1}$
associated with agent $y_{i}$'s proximity to $S$. 
We note that when $\mbox{int}\left(S\bigcap\mathcal{X}\right)$ is the empty set, 
$\phi_{i}(y_{i})$ is convex. 

\begin{figure}[h]
\begin{centering}
\includegraphics[width=0.3\columnwidth]{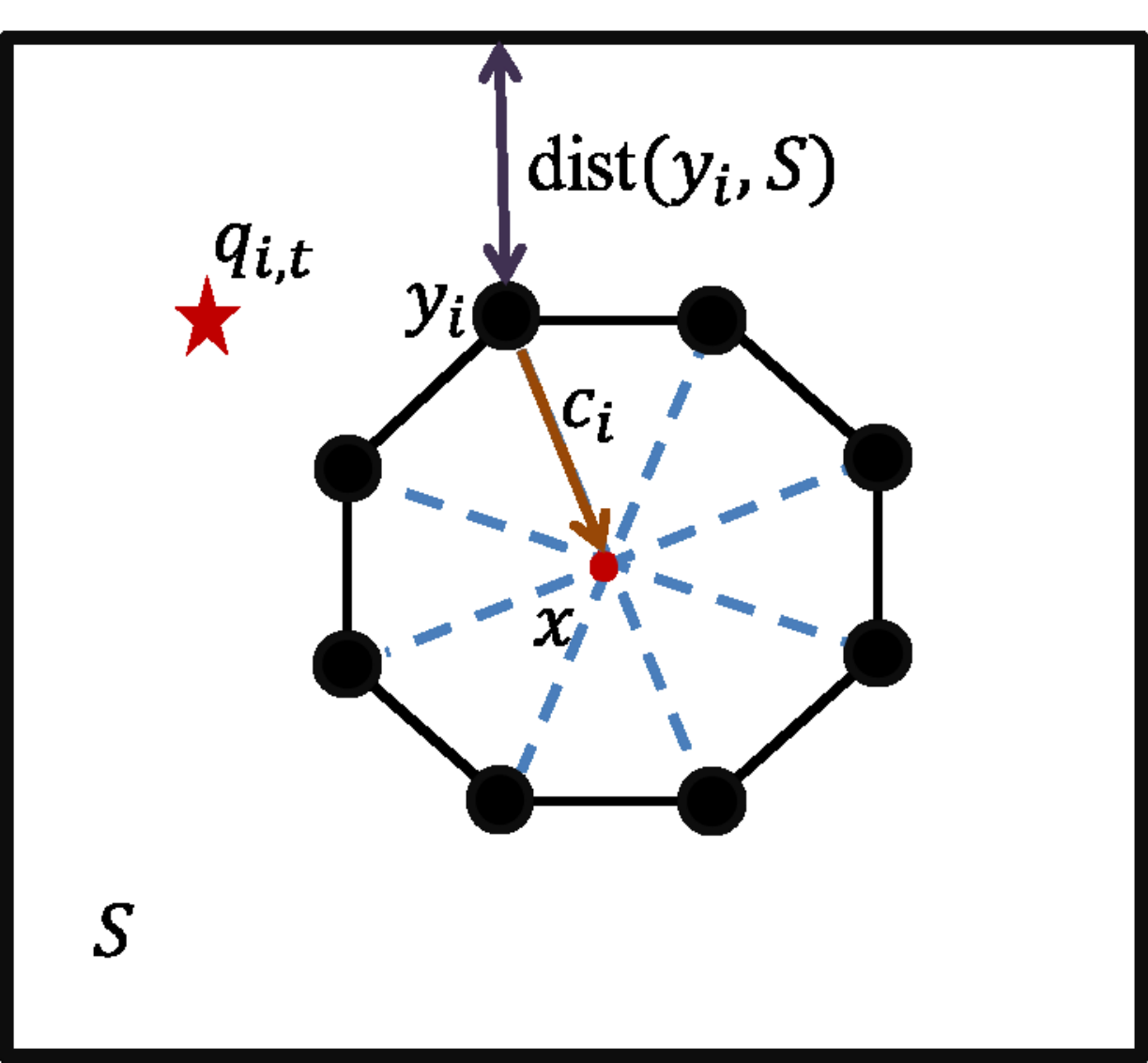}
\par\end{centering}

\protect\caption{\label{fig:Formation Acquisition Scenario}Formation acquisition problem
amongst six planar robots or agents.}
\end{figure}

At each time step $t$, agent $i$ obtains a location
of interest $q_{i,t}$ and the centroid is ideally located close to
these locations of interest promoted through the minimization of the
function $f_{i,t}(x)=\frac{1}{2}\left\Vert x-q_{i,t}\right\Vert _{2}^{2}$.
The example illustrated in Figure \ref{fig:Formation Acquisition Scenario}
takes the form of problem \eqref{eq:objective fun}, namely 
\begin{eqnarray*}
\min_{x\in\mathcal{X},y_{1},\dots y_{n}\in Y} & \sum_{t=1}^{T}\sum_{i=1}^{n}(f_{i,t}(x)+\phi_{i}(y_{i}))\\
\mbox{s.t.} & A_{i}x+B_{i}y_{i}=c_{i}\,\mbox{for all }i\in\left[n\right],
\end{eqnarray*}
where $A_{i}=-B_{i}=I_{2}$ for all $i\in\left[n\right]$.

Consider $S=\left\{ \left(x,y\right)\in\mathbb{R}|\left|x\right|=1.5,\left|y\right|=1.5\right\} $
and so $\phi_{i}(y_{i})=\left(2.5-\left\Vert y_{i}\right\Vert _{\infty}\right)^{-1}$.
The relevant parameters of the ADMM algorithm are $g_{i}(t)=\nabla f_{i,t}(x_{i})=x_{i}-q_{i,t}$,
$k=2,$ $\rho=0.5$, and $\psi(x)=\left\Vert x\right\Vert _{2}^{2}$.
The remaining terms of the regret bound are $L_{\phi}=4/9$, $L_{f}=\sqrt{2}$,
$\sigma_{1}(A_{i})=\sigma_{m_{i}}(B_{i})=1$, $D_{\lambda}=2$, and
$K=1$.

\begin{figure}[tb]
\begin{center}
\begin{subfigure}[b]{0.5\textwidth}
\includegraphics[scale=0.5]{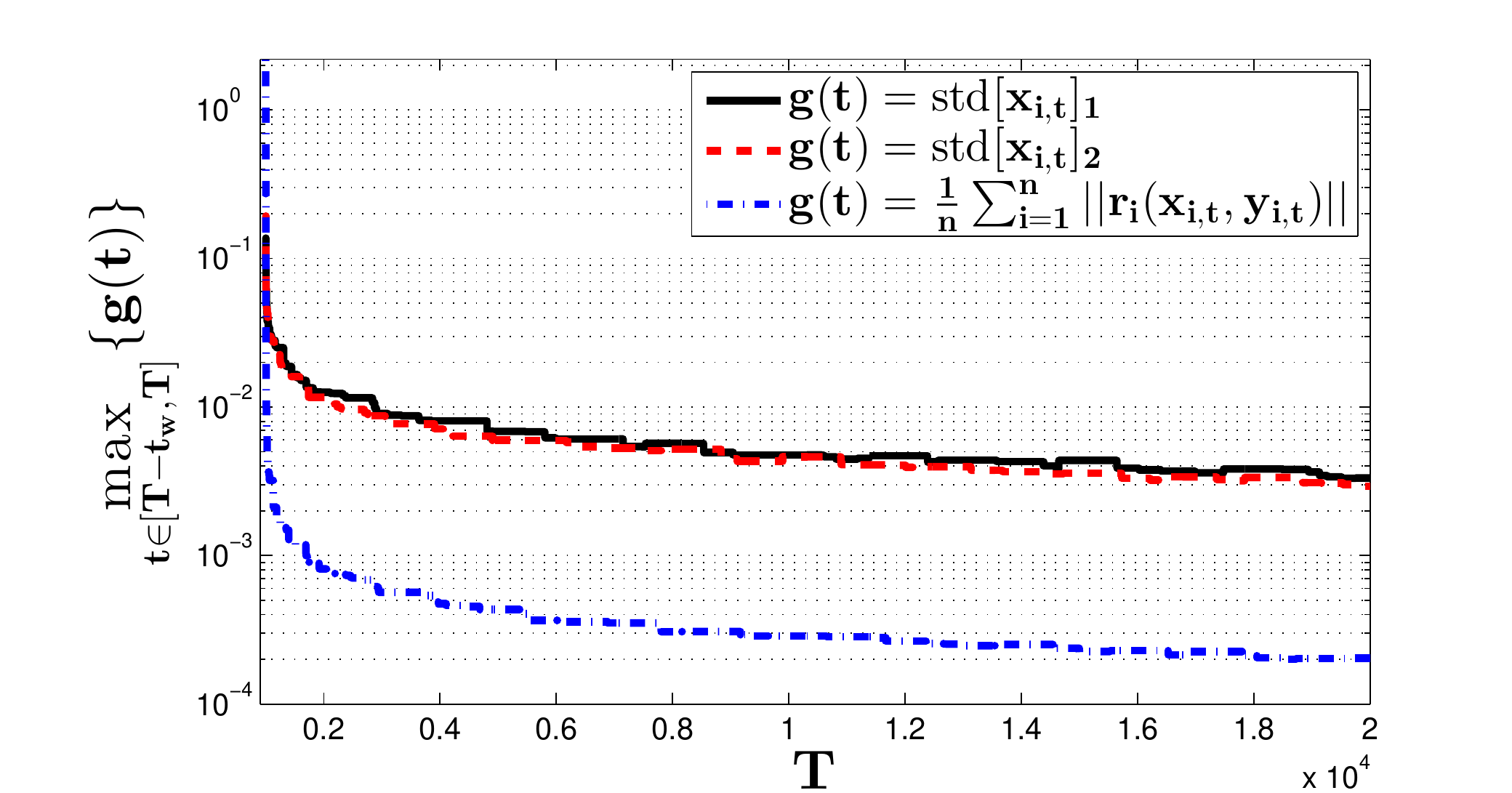}
\caption{Distributed Dual Averaging}
\includegraphics[scale=0.5]{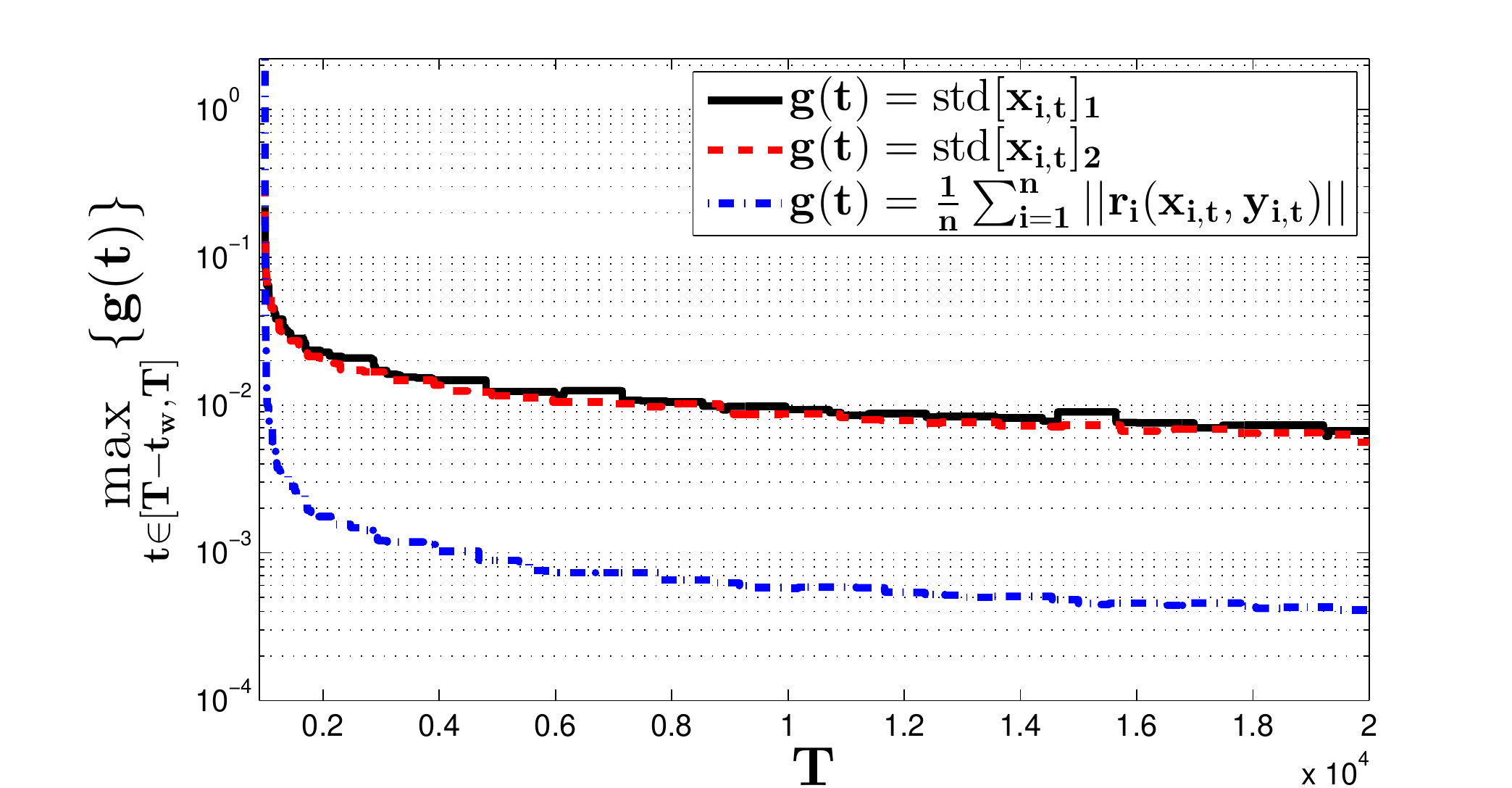}
\caption{Distributed Gradient Descent}
\end{subfigure}
\end{center}
\caption{\label{fig:Regret}The standard deviation of the global variable $x_{i}$
and the average residue for each agent over times smoothed by taking
the maximum over a $t_{w}=1000$ sliding window.}
\end{figure}

\begin{figure}[tb]
\begin{center}
\includegraphics[width=0.5\columnwidth]{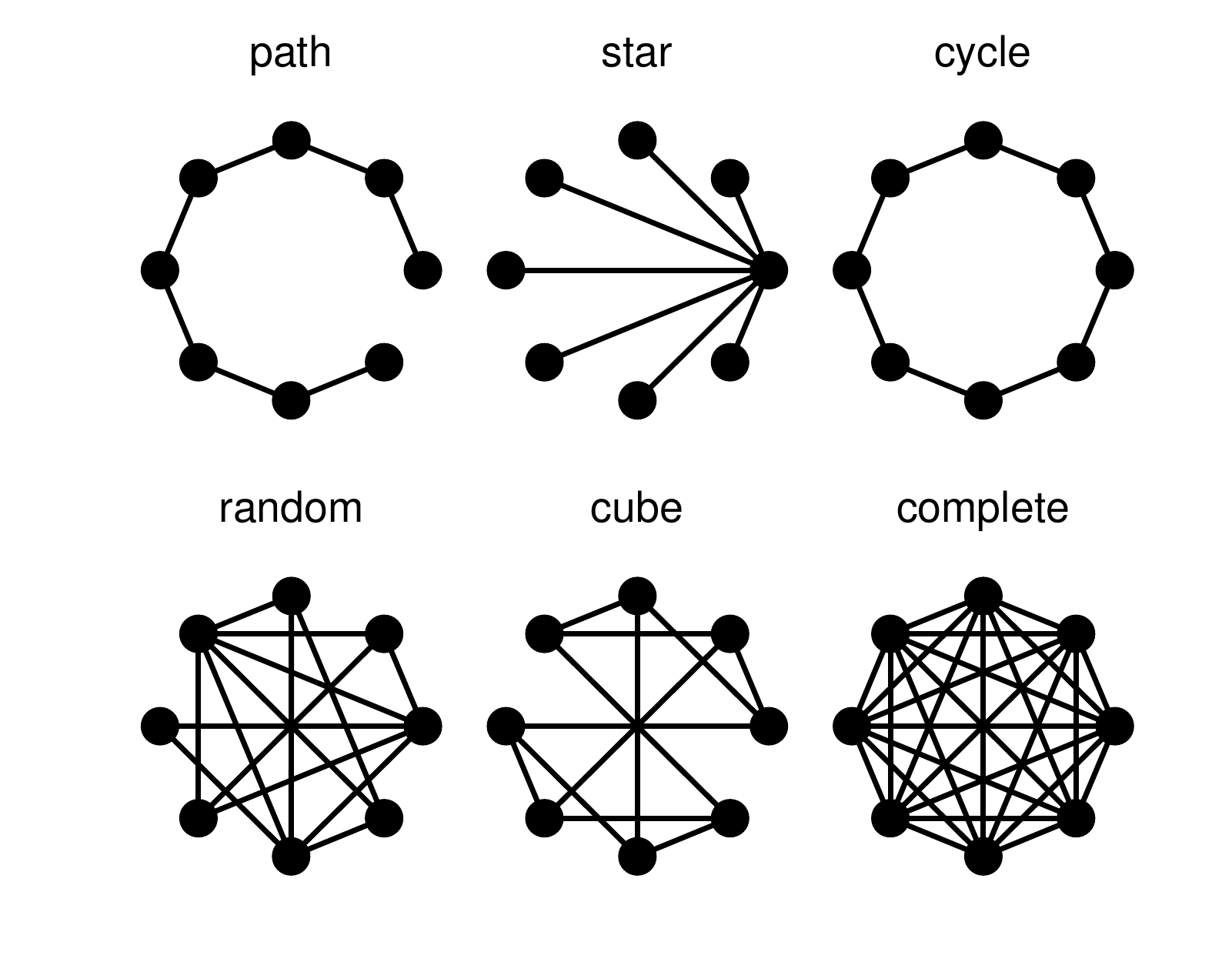}
\end{center}
\caption{\label{fig:GraphTopologies}Topologies of the six different graph types considered for the formation acquisition problem.}
\end{figure}

The algorithm was applied to $n=8$ agents connected over a random
graph (see Figure \ref{fig:GraphTopologies}) with $\sigma_{2}(P)=0.78$
with $c_{i}$'s selected to acquire a formation with $n$ agents equidistant
apart on the circumference of a circle of radius 0.4. Locations of
interest switch at each time step between a uniform distribution over
the area of a length 0.5 square centered at $\left(-0.75,0\right)$
and a Gaussian distribution with mean $\left(0,-0.75\right)$ and
standard deviation $0.01I_{2}$, with bounds outside of $\mathcal{X}$ ignored. 
The convergence of the global variables $x_{i,t}$ to agreement as
well as the reduction of the residue over time are displayed in Figure \ref{fig:Regret}. 
Note that the local copies of the global variable converge faster to consensus 
using the distributed DA as compared with embedding the distributed GD for OD-ADMM.

\begin{figure}[tb]
\begin{center}
\includegraphics[bb=40bp 0bp 430bp 288bp,clip,width=0.5\columnwidth]{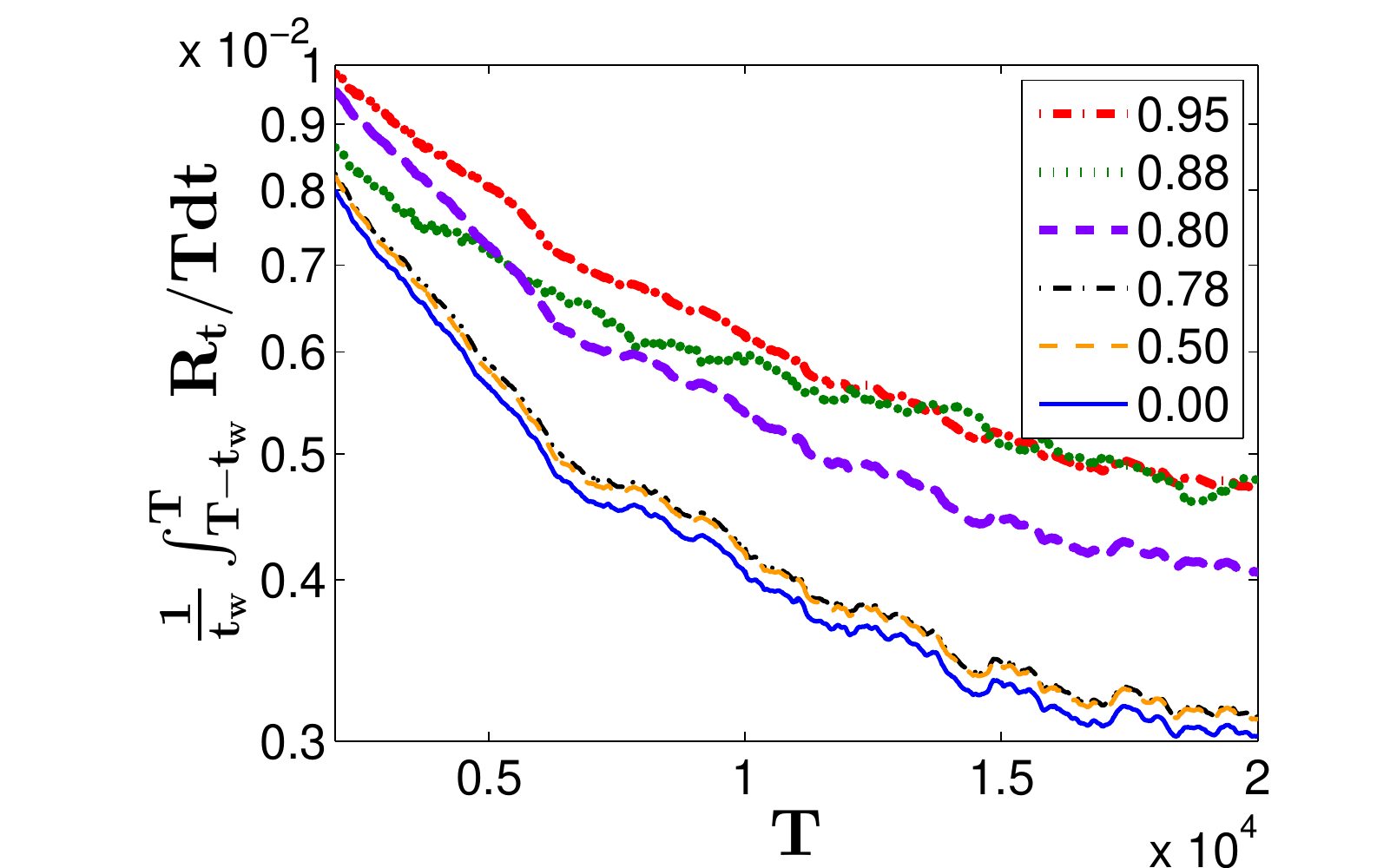}
\end{center}
\caption{\label{fig:GraphTypeComparison}The social regret per time $R_{T}/T$ performance
of six different graph types, specifically a path, star, cycle, random,
cube and complete graph with $\sigma_{2}(P)$ of $0.95,0.88,0.80,0.78,0.50,0.00$,
respectively. The trajectories are smoothed by taking the average
over a $t_{w}=1000$ sliding window. }
\end{figure}

The performance of the algorithm was compared for different graph
topologies, namely path, star, cycle, random, cube and complete graphs.
These graph topologies are displayed in Figure \ref{fig:GraphTopologies}.
The matrix $P$ was formed as proposed in Proposition \ref{prop:Strong},
with $\epsilon=d_{\max}+1$, as such $\sigma_{2}(P)=1-\frac{1}{\epsilon}\Lambda_{2}(L(\mathcal{G})).$
Under the same locations of interest as described previously, the performance
of the regret per time $R_{T}/T$ for each graph topology is compared
in Figure \ref{fig:GraphTypeComparison}. The performance strongly
correlates with $\sigma_{2}(P)$, as predicted in Theorem \ref{thm:Regret(xi)_GD},
with smaller $\sigma_{2}(P)$ exhibiting improved performance.

%
%
%
%
%
%
%
%
%
%
%
%

\section{Conclusion\label{sec:Conclusion}}

In this work, online distributed ADMM has been introduced and analyzed,
where a network of decision-makers or agents cooperatively
optimize an objective that decomposes to global and local objectives,
and is partially online.
Moreover, the local variables and the global variable are linearly constraint
(specific to each agent).
This problem setup has a wide range of applications in networked systems,
such as in distributed robotics and computer networks. A distributed
algorithm allows us
to make decisions across the network based on local data 
and information exchange with neighboring agents. 

The online distributed algorithm developed in this paper,
achieves a sub-linear social regret of $O(\sqrt{T})$, that simultaneously
captures sub-optimality of the objective function and the
violations in the linear local constraints.
In particular, this algorithm is competitive with respect to the best fixed decision performance
in hindsight. 
Moreover, we have highlighted the role of the underlying network
topology in achieving a ``good'' social regret, i.e., the regret bound improves
with increased connectivity in the network. 
The proposed algorithm was then applied to a formation acquisition problem.

Future work of particular interest includes exploring social regret over a time varying network,
and investigating favorable network characteristics for the proposed
online distributed ADMM algorithm.
\section{Appendix\label{sec:Appendix}}

The following results can be found in 
\cite{Duchi2012,Hosseini2013,Bertsekas2011};
as such they are presented here with no or abridged proofs.
\begin{prop}
\label{prop:Strong}If graph $\mathcal{G}$ is strongly connected
then the matrix $P=I-\frac{1}{\epsilon}\mbox{diag}\left(v\right)L(\mathcal{G})$
is doubly stochastic, where $v^{T}L\left(\mathcal{G}\right)=0$ with
positive vector $v=\left[v_{1},v_{2},\dots,v_{n}\right]^{T}$ and
$\epsilon\in(\max_{i\in V}(v_{i}d_{i}),\infty)$. If graph $\mathcal{G}$
is balanced then the matrix $P=I-\frac{1}{\epsilon}L(\mathcal{G})$
is doubly stochastic, where $\epsilon\in(d_{\max},\mbox{\ensuremath{\infty}})$. 
\end{prop}

\begin{prop}
\label{prop:orthogonal proj}For any $u\in\mathbb{R}^{m}$, $v\in\chi$,
and orthogonal projection operator $\prod_{\chi}$ onto $\chi$ we
have 
\[
\left\langle u-\prod_{\chi}(u),u-v\right\rangle \geq0.
\]

\end{prop}

\begin{lem}
\label{lem:x_i-y-1}For any $u,v\in\mathbb{R}^{m}$, and under the
conditions stated for the proximal function $\psi$ and step size $\alpha,$
we have 
\[
\Vert\prod_{\chi}^{\psi}(u,\alpha)-\prod_{\chi}^{\psi}(v,\alpha)\Vert\leq\alpha\Vert u-v\Vert_{*}.
\]

\end{lem}

\begin{lem}
\label{lem:DA network effect}For sequences $z_{i,t}$ and
$x_{i,t}$ generated by Algorithm\textup{ \ref{alg:OD-ADMM}}, using
the distributed DA method where,
\[
z_{i,t+1}=\sum_{j=1}^{n}P_{ji}z_{j,t}+g_{i,t}+A_{i}^{T}\lambda_{i,t+1}
\]
and $x_{i,t+1}=\prod_{\chi}^{\psi}(z_{i,t+1},\alpha_{t})$, we have
\begin{align*}
\Vert\theta_{t}-x_{i,t}\Vert_{*} & \leq\alpha_{t-1}\frac{\sqrt{n}(L_{f}+\Km)}{1-\sigma_{2}(P)}
\end{align*}
 for all $i\in[n]$ and $t\in[T]$, where the sequence $\theta_{t}$
is generated by \eqref{eq:theta update_DA}, $\K_{i}=L_{\phi} {\sigma_{1}(A_{i})/\sigma_{1}(B_{i})}$,
and $\Km=\max_{i}\mathcal{\K}_{i}$.\end{lem}
\begin{IEEEproof}
Based on the definition of $z_{i,t}$ we have
\begin{align*}
z_{i,t} & =\sum_{j=1}^{n}\left[P^{t-1}\right]_{ji}z_{j,1}+\sum_{k=1}^{t-1}\sum_{j=1}^{n}\left[P^{k-1}\right]_{ji}\left(g_{j,t-k}+A_{j}^{T}\lambda_{j,t-k+1}\right).
\end{align*}
In addition, $z_{t}$ evolves as
\begin{align}
z_{t} & =z_{1}+\sum_{k=1}^{t-1}\sum_{j=1}^{n}\frac{1}{n}(g_{j,t-k}+A_{j}^{T}\lambda_{j,t-k+1}).\label{eq:q_bar}
\end{align}
Assuming $z_{i,1}=0$ for all $i\in[n]$ and based on \eqref{eq:q_bar}
we have
\begin{align}
z_{t}-z_{i,t} & =\sum_{k=1}^{t-1}\sum_{j=1}^{n}(\frac{1}{n}-\left[P^{k-1}\right]_{ji})(g_{j,t-k}+A_{j}^{T}\lambda_{j,t-k+1})\label{eq:z-z-i-evol}
\end{align}
Thus, the dual norm of $z_{t}-z_{i,t}$ can be bounded as 
\begin{align}
 \Vert z_{t}-z_{i,t}\Vert_{*} & \leq\sum_{k=1}^{t-1}\sum_{j=1}^{n}\left\Vert g_{j,t-k}+A_{j}^{T}\lambda_{j,t-k+1}\right\Vert _{*}|\frac{1}{n}-\left[P^{k-1}\right]_{ji}|\nonumber \\
 & \leq\sum_{k=1}^{t-1}\max_{j}\left\Vert g_{j,t-k}+A_{j}^{T}\lambda_{j,t-k+1}\right\Vert _{*}\Vert P^{k-1}e_{i}-\frac{\mathbf{1}}{n}\Vert_{1}.\label{eq:dual (q-q_i evol)}
\end{align}
Since $\left\Vert g_{i,t}\right\Vert _{*}\leq L_{f}$ and $||A_{i}^{T}\lambda_{i,t}||_{*}\leq\mathcal{\K}_{i}\leq\Km$,
the dual norm of $z_{t}-z_{i,t}$ is further bounded as %
\footnote{Note that $\Vert P^{t}x-\frac{\mathbf{1}}{n}\Vert_{1}\leq\sigma_{2}(P)^{t}\sqrt{n},$
where the vector $x$ belongs to $\{x\in\mathbb{R}^{n}\vert 
x\geq0,\sum_{i=1}^{n}x_{i}=1\}$; this property of stochastic matrices was similarly used by Duchi \textit{et
al.} \cite{Duchi2012}. %
}
\begin{equation}
\Vert z_{t}-z_{i,t}\Vert_{*}\leq\sqrt{n}(L_{f}+\Km)\sum_{k=1}^{t-1}\sigma_{2}(P)^{k-1}.\label{eq: q-q_i-2}
\end{equation}
In addition, as $P$ is a doubly stochastic matrix, $\sigma_{2}(P)\leq1$~\cite{Berman1979}. Thus, the
inequality \eqref{eq: q-q_i-2} is bounded as
\begin{align*}
\Vert z_{t}-z_{i,t}\Vert_{*} & \leq\frac{\sqrt{n}(L_{f}+\Km)}{1-\sigma_{2}(P)},
\end{align*}
Since $\theta_{t}=\prod_{\chi}^{\psi}(z_{t},\alpha_{t-1})$ and $x_{i,t}=\prod_{\chi}^{\psi}(z_{i,t},\alpha_{t-1})$,
the statement of the lemma follows from Lemma \ref{lem:x_i-y-1}.\end{IEEEproof}
\begin{lem}
\label{lem:GD network effect}For sequences $x_{i,t}$ and
$h{}_{i,t}$ generated by Algorithm\textup{ \ref{alg:OD-ADMM} using
distributed GD method,} where 
\[
h_{i,t}=\sum_{j=1}^{n}P_{ji}x_{j,t-1}-\alpha_{t-1}(g_{i,t-1}+A_{i}^{T}\lambda_{i,t})
\]
and $x_{i,t}=\prod_{\chi}h_{i,t}$, we have 
\begin{align*}
\Vert\theta_{t}-x_{i,t}\Vert\leq2\sqrt{n}(L_{f}+\Km)\sum_{k=1}^{t-1}\alpha_{t-k}\sigma_{2}(P)^{k-1}
\end{align*}
 for all $i\in[n]$ and $t\in[T]$, where the sequence $\theta_{t}$
is generated by \[
\theta_{t}=\frac{1}{n}\sum_{i=1}^{n}x_{i,t},
\]
$\K_{i}=L_{\phi}{\sigma_{1}(A_{i})}/{\sigma_{1}(B_{i})}$, and $\Km=\max_{i}\mathcal{\K}_{i}$.\end{lem}
\begin{IEEEproof}
Denote $r_{i,t}=x_{i,t}-h_{i,t}$; thus based on the definition
of $h_{i,t}$ we have
\begin{align}
x_{i,t} =h_{i,t}+r_{i,t}  =\sum_{j=1}^{n}P_{ji}x_{j,t-1}-\alpha_{t-1}(g_{i,t-1}+A_{i}^{T}\lambda_{i,t})+r_{i,t}\label{eq:primal_up_Gd}
\end{align}
 Subsequently, we can represent $x_{i,t}$ as 
\begin{align*}
x_{i,t} & =\sum_{j=1}^{n}\left[P^{t-1}\right]_{ji}x_{j,1}-\sum_{k=1}^{t-1}\sum_{j=1}^{n}\left[P^{k-1}\right]_{ji}(\alpha_{t-k}\left(g_{j,t-k}+A_{j}^{T}\lambda_{j,t-k+1}\right)-r_{j,t-k+1}).
\end{align*}
 In addition, based on \eqref{eq:primal_up_Gd}, the average primal
variable $\theta_{t}$ evolves as
\begin{align}
\theta_{t} & =\theta_{1}-\sum_{k=1}^{t-1}\sum_{j=1}^{n}\frac{1}{n}(\alpha_{t-k}(g_{j,t-k}+A_{j}^{T}\lambda_{j,t-k+1})-r_{j,t-k+1}).\label{eq:q_bar-1}
\end{align}
Assuming $x_{i,1}=0$ for all $i\in[n]$ and based on \eqref{eq:q_bar-1},
we can represent the network effect, that is the difference between
the average primal variable over the network and individual primal
variables, as 
\begin{align}
\theta_{t}-x_{i,t}= & \sum_{k=1}^{t-1}\sum_{j=1}^{n}(\frac{1}{n}-\left[P^{k-1}\right]_{ji})(r_{j,t-k+1}-\alpha_{t-k}(g_{j,t-k}+A_{j}^{T}\lambda_{j,t-k+1})).\label{DUPLICATE: eq:z-z-i-evol-1}
\end{align}
 Thus, the network effect \eqref{DUPLICATE: eq:z-z-i-evol-1} can
be bounded as 
\begin{align}
\Vert\theta_{t}-x_{i,t}\Vert \leq & \sum_{k=1}^{t-1}\sum_{j=1}^{n}\left\Vert r_{j,t-k+1}-\alpha_{t-k}(g_{j,t-k}+A_{j}^{T}\lambda_{j,t-k+1})\right\Vert |\frac{1}{n}-\left[P^{k-1}\right]_{ji}|\nonumber \\
\leq & \sum_{k=1}^{t-2}\max_{j}\left\Vert r_{j,t-k+1}-\alpha_{t-k}(g_{j,t-k}+A_{j}^{T}\lambda_{j,t-k+1})\right\Vert \Vert P^{k-1}e_{i}-\frac{1}{n} \mathbf{1}\Vert_{1}\label{eq:dual (q-q_i evol)-1}
\end{align}

Moreover, the difference between $h_{i,t}$ and its projection onto
$\chi$ is bounded as 
\begin{align*}
||r_{i,t}|| & =||\prod_{\chi}h_{i,t} -h_{i,t}||\leq||\sum_{j=1}^{n}P_{ji}x_{j,t-1}-h_{i,t}||\\
 & \leq\alpha_{t-1}||g_{i,t-1}+A_{i}^{T}\lambda_{i,t}||.
\end{align*}
Since $\left\Vert g_{i,t}\right\Vert _{*}\leq L_{f}$ and $||A_{i}^{T}\lambda_{i,t}||\leq\mathcal{\K}_{i}\leq\Km$,
the network effect is further bounded as %
\begin{equation}
\Vert\theta_{t}-x_{i,t}\Vert\leq2\sqrt{n}(L_{f}+\Km)\sum_{k=1}^{t-1}\alpha_{t-k}\sigma_{2}(P)^{k-1} . \label{eq: q-q_i-2-1}
\end{equation}

\end{IEEEproof}

\begin{lem}
\label{lem:DA sub-optimality effect}For any positive and non-increasing
sequence $\alpha(t)$ and $x^{*}\in\chi$
\begin{align*}
\sum_{t=1}^{T}\left\langle g_{t},\theta_{t}-x^{*}\right\rangle  & \leq\sum_{t=1}^{T}\frac{1}{n}\sum_{i=1}^{n}\left\langle A_{i}^{T}\lambda_{i,t+1},x^{*}-\theta_{t}\right\rangle +\frac{1}{\alpha_{T}}\psi(x^{*})+\sum_{t=2}^{T}\frac{\alpha_{t-1}}{2} \, ||g_{t}+\frac{1}{n}\sum_{i=1}^{n}A_{i}^{T}\lambda_{i,t+1}||_{*}^{2},
\end{align*}
 where the sequence $\theta_{t}$ is generated by \eqref{eq:z_bar update_DA}-\eqref{eq:theta update_DA}.\end{lem}
\begin{IEEEproof}
Based on Lemma 3 in \cite{Duchi2012}, we have 
\begin{align*}
\sum_{t=1}^{T}\left\langle g_{t}+\frac{1}{n}\sum_{i=1}^{n}A_{i}^{T}\lambda_{i,t+1},\theta_{t}-x^{*}\right\rangle  & \leq\frac{1}{\alpha_{T}}\psi(x^{*})+\sum_{t=2}^{T}\frac{\alpha_{t-1}}{2} \, ||g_{t}+\frac{1}{n}\sum_{i=1}^{n}A_{i}^{T}\lambda_{i,t+1}||_{*}^{2},
\end{align*}
and the statement of the lemma follows. 
\end{IEEEproof}

\begin{lem}
\label{lem:GD sub-optimality effect}For any positive and non-increasing
sequence $\alpha(t)$ and $x^{*}\in\chi$
\begin{align*}
\sum_{t=1}^{T}\left\langle g_{t},\theta_{t}-x^{*}\right\rangle  & \leq\frac{1}{n}\sum_{t=1}^{T}\sum_{j=1}^{n}\left\langle A_{j}^{T}\lambda_{j,t+1},x^{*}-\theta_{t}\right\rangle +\frac{2}{n^{2}}\sum_{t=1}^{T}\alpha_{t}(\sum_{j=1}^{n}||g_{j,t}+A_{j}^{T}\lambda_{j,t+1}||)^{2}\\
 & + \frac{1}{2\alpha_{1}}D_{\chi}^{2}+\sum_{t=1}^{T}(4L_{f}\sum_{k=0}^{t-1}\alpha_{t-k}\sigma_{2}(P)^{k})\frac{1}{n}\sum_{j=1}^{n}||(g_{j,t}+A_{j}^{T}\lambda_{j,t+1}||,
\end{align*}
 where the sequence of $\theta_{t}$ is generated by \eqref{eq:GD local h update},
\eqref{eq:GD local x update}, and \eqref{eq:theta update_GD}.\end{lem}
\begin{IEEEproof}
Denote $r_{i,t}=x_{i,t}-h_{i,t}$; thus based on the definition
of $h_{i,t}$, we have
\begin{align}
x_{i,t+1} & =h_{i,t+1}+r_{i,t+1}\nonumber \\
 & =\sum_{j=1}^{n}P_{ji}x_{j,t}-\alpha_{t}(g_{i,t}+A_{i}^{T}\lambda_{i,t+1})+r_{i,t+1} . \label{eq:primal_up_Gd-1}
\end{align}
 Subsequently, based on \eqref{eq:theta update_GD}, the average primal
variable $\theta_{t}$ evolves as
\begin{align}
\theta_{t+1} & =\theta_{t}-\frac{1}{n}\sum_{j=1}^{n}(\alpha_{t}(g_{j,t}+A_{j}^{T}\lambda_{j,t+1})-r_{j,t+1}).\label{eq: GD theta evolve}
\end{align}
Now, we can represent the deviation of average primal variable $\theta_{t}$
from $x^{*}$ as
\begin{align}
||\theta_{t+1}-x^{*}||^{2}= & ||\theta_{t}-x^{*}||^{2}+\frac{1}{n^{2}}||\sum_{j=1}^{n}(\alpha_{t}(g_{j,t}+A_{j}^{T}\lambda_{j,t+1})-r_{j,t+1})||^{2}\nonumber \\
 & - \frac{2\alpha_{t}}{n}\sum_{j=1}^{n}\left\langle (g_{j,t}+A_{j}^{T}\lambda_{j,t+1}),\theta_{t}-x^{*}\right\rangle +\frac{2}{n}\sum_{j=1}^{n}\left\langle r_{j,t+1},\theta_{t}-x^{*}\right\rangle .\label{DUPLICATE: eq:theta-x^*}
\end{align}
 
Note that 
\begin{align*}
||\sum_{j=1}^{n}(\alpha_{t}(g_{j,t}+A_{j}^{T}\lambda_{j,t+1})-r_{j,t+1})||^{2} & \leq \left\{ \sum_{j=1}^{n}||\alpha_{t}(g_{j,t}+A_{j}^{T}\lambda_{j,t+1}||+||r_{j,t+1}||\right\}^{2}\\
 & \leq4\alpha_{t}^{2}(\sum_{j=1}^{n}||g_{j,t}+A_{j}^{T}\lambda_{j,t+1}||)^{2},
\end{align*}
and from Proposition \ref{prop:orthogonal proj} we have 
\begin{align*}
\left\langle r_{j,t+1},\theta_{t}-x^{*}\right\rangle  & \leq\left\langle r_{j,t+1},\theta_{t}-h_{j,t+1}\right\rangle +\left\langle \prod_{\chi}h_{j,t+1}-h_{j,t+1},h_{j,t+1}-x^{*}\right\rangle \\
 & \leq\left\langle r_{j,t+1},\theta_{t}-h_{j,t+1}\right\rangle \leq\alpha_{t}||(g_{j,t}+A_{j}^{T}\lambda_{j,t+1}||||\theta_{t}-h_{j,t+1}||.
\end{align*}
Based on Lemma 8 in \cite{Yan2010}, we have 
\begin{align*}
\left\langle r_{j,t+1},\theta_{t}-x^{*}\right\rangle  & \leq4L_{f} \alpha_{t}||(g_{j,t}+A_{j}^{T}\lambda_{j,t+1}||\sum_{k=0}^{t-1}\alpha_{t-k}\sigma_{2}(P)^{k}.
\end{align*}
Thus, by rearranging the terms in \eqref{DUPLICATE: eq:theta-x^*},
we have 
\begin{align}
\sum_{t=1}^{T}\left\langle g_{t},\theta_{t}-x^{*}\right\rangle \leq & \frac{1}{2\alpha_{1}}||\theta_{1}-x^{*}||^{2}-\frac{1}{2\alpha_{T+1}}||\theta_{T+1}-x^{*}||^{2}+\frac{2}{n^{2}}\sum_{t=1}^{T}\alpha_{t}(\sum_{j=1}^{n}||g_{j,t}+A_{j}^{T}\lambda_{j,t+1}||)^{2}\nonumber \\
 & + \frac{1}{n}\sum_{t=1}^{T}\sum_{j=1}^{n}\left\langle A_{j}^{T}\lambda_{j,t+1},x^{*}-\theta_{t}\right\rangle +\sum_{t=1}^{T}(4L_{f}\sum_{k=0}^{t-1}\alpha_{t-k}\sigma_{2}(P)^{k})\frac{1}{n}\sum_{j=1}^{n}||(g_{j,t}+A_{j}^{T}\lambda_{j,t+1}||.\label{DUPLICATE: eq:theta-x^*-1}
\end{align}

Since the diameter of $\chi$ is bounded by $D_{\chi}$, we have $||\theta_{1}-x^{*}||^{2}\leq D_{\chi}^{2}$
and the statement of the lemma follows.
\end{IEEEproof}
\bibliographystyle{IEEEtran}
\bibliography{DSSL-OnlineRegret-ADMM,DSSL-OnlineRegret-CDC2013,generalPapers,regretPapers,cdc2013_dist_opt,DSSL-OnlineRegret,ADMM_journal2014}
\end{document}